\documentclass[11pt]{article}
\usepackage[utf8]{inputenc}

\title{A multi-objective optimization framework for sustainable transitions}

\author{Cris R. Hasan$^{1*}$ \and Luigi Cao Pinna$^2$ \and John Crawford$^1$ \and Stuart Kauffman$^3$ \and Roger Koppl$^4$ \and Jonathan Lee$^1$ \and Demival Vasques$^5$ \and Edward Weinberger$^6$}

\date{%
    $^1$University of Glasgow, UK, $^2$Research Centre for Agricultural and Environment, Italy, $^3$Institute for Biocomplexity and Informatics, US, $^4$Syracuse University, US, \\$^5$University of Luxembourg, Luxembourg, $^6$New York University, US \\[2ex]%
    $^*$Corresponding author. Email: cris.hasan@glasgow.ac.uk.
}

\usepackage[utf8]{inputenc}

\usepackage{natbib}
\usepackage{graphicx}
\usepackage{amssymb}
\usepackage{amsmath}
\usepackage{hyperref}
\usepackage{comment}
\usepackage[capitalise]{cleveref}
\usepackage{amsthm}
\usepackage{subcaption}
\usepackage[font=small,labelfont=bf]{caption}
\crefformat{equation}{(#2#1#3)}
\usepackage[dvipsnames]{xcolor}
\usepackage{todonotes}
\usepackage[left=2cm]{geometry}
\usepackage[normalem]{ulem}
\usepackage{makecell}
\definecolor{mygray}{gray}{0.6}
\DeclareMathOperator*{\argmax}{arg\,max}
\usepackage{appendix}

\hypersetup{
    colorlinks=true,
    citecolor=blue,      
    linkcolor=red,       
    urlcolor=magenta     
}

\usepackage{caption}
\begin{document}

\maketitle

\section*{Abstract}
Achieving a just and sustainable transition requires the pursuit of multiple social and environmental targets.
Two primary barriers impede this process: (1) targets are often in conflict with each other, and (2) policies aimed at these targets are commonly planned in isolation, neglecting complex interdependencies in the system. 
To address these challenges, we propose a general modeling framework that evaluates the holistic impact of policies and decision-making on sustainability targets while capturing system interdependencies in a policy-target network.
Inspired by Kauffman's NK fitness landscape, our framework takes the form of a multi-objective optimization model that employs a dynamic evolutionary algorithm in conjunction with network analysis.
Our algorithm accounts for tradeoffs between conflicting targets by dynamically reallocating resources 
to the most impactful and efficient policies.
One key finding indicates that increasing resources generally enhances performance, but marginal gains stagnate at a point of diminishing returns.
Sensitivity analysis reveals that the system is primarily driven by three factors: budget constraint, network density (interconnectivity), and policy efficacy.
This study serves as a foundational step towards developing a decision-support tool that assists policymakers in achieving optimal outcomes for problems with a large number of dynamically interacting targets.

\vspace{0.5cm}

\noindent\textbf{Keywords:} Sustainable Development Goals, Multi-objective optimization, Fitness landscapes, Decision-making and policy

\pagebreak

\section{Introduction}\label{sec1}

A decade after the adoption of the Paris Agreement \citep{agreement2015paris}, climate change remains one of the most pressing challenges of the 21st century. 
Despite growing global momentum towards a net-zero transition, environmental and social policies often operate in silos \citep{buylova2025bridging}.
This compartmentalized planning increases the risk of policy misalignment and missed opportunities in cultivating 
mutually reinforcing outcomes \citep{lah2025breaking}.
What is more, sustainability targets can sometimes be in conflict with justice goals \citep{abram2022just, ciplet2020transition}. 
Therefore, it is crucial to ensure that future transitions are designed to be not only environmentally sustainable but also equitable and socially just.

In response to these challenges, Doughnut Economics \citep{raworth2017doughnut} offers an alternative economic theory that holistically integrates social and ecological challenges, by synthesizing \emph{Sustainable Development Goals} (SDGs) \citep{SDGs2015} with the \emph{Planetary Boundaries} framework \citep{rockstrom2009planetary}.
This theory highlights the urgent need for a paradigm shift in political and economic mindsets when pursuing an equitable system transformation.
Recent studies applied Doughnut Economics theory to articulate aspirational vision for a set of social and environmental goals in various urban contexts; see for example \cite{acosta2022linking, cattaneo2025ecological, hjelmskog2023thriving}. 
However, given the intricate complexities of interlinkages in the system, a full operationalisation of the this theoretical framework necessitates a whole-system optimization approach to dynamically reconcile tradeoffs among competing goals.

In this paper, we propose a conceptual model for optimizing social, economic and environmental outcomes through decision-making process. 
Drawing from the 169 targets outlined in the SDG framework, we adopt the term \emph{sustainability targets}, or simply \emph{targets}, to denote specific objectives that represent desired environmental and wellbeing outcomes.
We further define the term \emph{policies} to broadly encompass a suite of policy levers available to decision- and policy-makers in pursuit of sustainability targets \citep{chan2020levers}.
Examples of policy levers include, but are not limited to, fiscal measures (e.g., taxes and subsidies),   regulatory interventions, and infrastructure investments \citep{schlesier2024measuring}.
The dynamic interactions between policies and targets will be formalized as a two-mode network, which we term the \emph{policy-target network}.

To reconcile tradeoffs among multiple competing sustainability targets, we adopt a multi-objective optimization approach. 
Specifically, we apply a generalized NK fitness model approach \citep{kauffman1989nk}
in conjunction with an evolutionary algorithm \citep{deb2011multi} and network analysis to construct a conceptual model for the purpose of optimizing sustainability targets.
The NK fitness model \citep{kauffman1989nk, kauffman1992origins} was originally developed to study complex systems (e.g., biological systems) and explore how the overall performance of such systems (e.g., species fitness) can be optimized.
The model's two key parameters $N$ and $K$ represent the number of elements and their interconnectedness in the system, respectively.
Systems with few interactions (low $K$) are generally easier and faster to optimize but cannot capture the complexities of real-world interactions. 
Conversely, fitness landscapes with a larger number of interactions (high $K$) are more difficult to navigate but remain robust to perturbations.
Therefore, an intermediate value of $K$ is favorable when seeking balance between accelerating the transition and maintaining a degree of robustness and adaptability \citep{bull2022nonbinary, nowak2015analysis, srivastava2023alphabet}.

The conceptualization of our framework is inspired by the NK fitness landscape in three primary ways.
Firstly, in classical fitness landscapes, a genome is represented by a discrete string of genes.
Analogous to this concept, we introduce the term \emph{policy array}, defined as an ordered sequence of policy levers.
As opposed to the binary case of the canonical NK model, we employ a non-binary \emph{alphabet} cardinality \citep{bull2022nonbinary, srivastava2023alphabet} to accommodate multi-level resource deployments for policies. 
Secondly, the notion of species fitness serves as an objective function in the original model.
We borrow this notion to quantify prosperity \citep{jackson2009prosperity} or thriving \citep{raworth2017doughnut} and operationalize these concepts as a weighted sum of a diverse set of environmental and wellbeing outcomes.
The third essential component of the NK fitness landscape is the mutation process of genomes, which can be visualized as the navigation of a landscape of peaks and valleys that represent different fitness levels.
Akin to this process is a dynamic sequence of policy interventions by which the state of the system is incrementally adjusted in search of an optimal solution within a \emph{performance landscape}.

Despite recent advances \citep{stechemesser2024climate}, there is no standard practice for evaluating policies or measuring the impact of interventions on intended deliverables, let alone for capturing unforeseen unintended consequences.
To circumvent this challenge, we analyze our optimization model across a spectrum of hypothetical scenarios.
Specifically, we leverage \emph{configuration models} \citep{bollobas1980probabilistic}---a class of statistical techniques that enables the systematic generation and exploration of a broad range of network topologies.
These models will also be employed to assign weights for the network interconnections, thereby quantifying the influence of policy levers on the impacted targets. 
Furthermore, targets are ranked according to their relative importance by configuring appropriately specified random distributions. 
This hybrid approach enables our model to provide valuable insights to 
decision-makers, especially in contexts where comprehensive empirical data are unavailable.

By investigating the performance of the system under a diverse set of statistical scenarios,  
we show that increasing resources generally enhances system performance, but the rate of marginal improvements plateaus at higher budgets, indicating a point of diminishing returns.
This phenomenon has been reported in various relevant contexts \citep{raworth2017doughnut, jackson2009prosperity, stern2007economics}.
Furthermore, systems with fewer policies stagnate more rapidly with increasing budgets, while incorporating expanded policy frameworks demonstrated enhanced resilience through risk diversification and adaptive capacity.
Our study reveals that systems with high network interconnectivity demonstrate resource efficiency, yet their optimization necessitates an intricate consideration of inherent trade-offs among targets.
We also carry out simulation-based sensitivity analysis to examine the system's response to variations in input parameters; this comprehensive approach allows us to identify the key factors that govern the qualitative and quantitative behavior of the system.

\begin{table}[t!]
\centering
\begin{tabular}{|l|l|l|l|l}
 
\textbf{Parameter} & \textbf{Description} & \textbf{Baseline Value} & \textbf{Range} &  \\  
$N$ & Number of policies  & 100 &  $\mathbb{Z}^+$  &  \\  
$M$ & Number of targets  & 30 &  $\mathbb{Z}^+$   &  \\  
$A$ & Alphabet  & 5 &   $\mathbb{Z}^+/\{1\}$  &  \\  
$B$ & Average per-policy budget  & 3 &   $[0,A-1]$  &  \\  
$B_T$ & Total budget given by $N \times B$  & 300 &   $[0,N(A-1)]$  &  \\  
$\mu_k$ & Average outdegree of policies & 5 &   $(1,M]$  &  \\  
$\rho$ & Network density estimated by $\mu_k/M$  & 1/6 &   $(1/M,1]$  &  \\  
$\beta_k$ & Outdegree scale parameter  & 2 &   $\mathbb{R}^+$  &  \\  
$\mu_c$ & \makecell[l]{Average signed-coefficient (also \\ representative of the overall policy efficacy)}  & 1/3 &   $\mathbb{R}^+$  &  \\  
$\beta_c$ & Signed-coefficient scale parameter  & 2 &   $\mathbb{R}^+$  &  \\  
$\mu_w$ & Relative-importance average parameter  & 8 &   $\mathbb{R}^+$  &  \\  
$\beta_w$ & Relative-importance scale parameter  & 15 &   $\mathbb{R}^+$  &  \\  
$\eta$ & Scale parameter for $\mathcal{H}$  & 3 &   $\mathbb{R}^+$  &  \\  
\end{tabular}
\caption{Table of system parameters with their values and range.}
\label{table:parameters}
\end{table}

\section{Methodology}
\label{sec:Methodology}
Inspired by Kauffman's fitness landscape model \citep{kauffman1989nk, kauffman1992origins}, our analytic framework combines multi-objective combinatorial optimization, evolutionary algorithm \citep{deb2011multi}, and network science. 
In classical NK fitness landscapes, the fixed prescribed parameters $N$ and $K$ denote the number of genes and the degree of their interdependence, respectively. 
In our adaptation, $N$ signifies the number of policies or decisions.
We further extend the model by introducing a new parameter $M$, which denotes the number of targets associated with desired social and environmental outcomes. 
The stochastic quantity $K_i^{out}$ captures the outdegree of a policy $i$, $\forall i \in {1,2,...,N}$.
The parameter $\mu_k$ represents the average number of targets influenced by a single policy.
A summary of the primary model parameters is provided in    \cref{table:parameters}.

In contrast to the classical NK model which assumes that the building blocks of the model (e.g., genes) are binary, we adopt a more nuanced perspective. 
Namely, we allow for a multi-level resource allocations by fine tuning the value of each policy within a broader discrete non-binary alphabet denoted as $A$ \citep{bull2022nonbinary, srivastava2023alphabet}, which represents the cardinality of the permissible numerical range for each policy.

To account for resource constraints, we introduce a budget parameter $B$, which is the average of resource deployment per policy accommodated by the total available budget $B_{T} = B N$.
Furthermore, analogous to the fitness function in traditional landscape models, we define a set of \emph{individual performance functions} that evaluate the progress of each target.
These functions are subsequently aggregated via a weighted function that quantifies the overall performance of the system.

In the following subsections, we first introduce the structure of the policy-target interactions, which comes in a form of a \emph{directed two-mode network}, followed by an overview of our approach for constructing this network using random distributions. 
We then describe the allocation of resources to an array of policies within a finite budget.
This is followed by a description of our method of computing the individual performance function for each target as well as the overall performance while considering the hierarchical relative importance among targets.
Then, we detail the evolutionary algorithm for our multi-objective optimization.
The section is concluded with a heuristic rationale for the selection of parameter values. 
For a quick familiarization of the proposed framework, please see the glossary listed in \cref{secA1}. 

\subsection{A directed two-mode network structure}
Two-mode networks are a specific type of network in which nodes are categorized into two distinct sets.
Connections, also known as links, or edges, exist strictly between nodes belonging to different sets.
Mathematically, a two-mode network is defined as a \emph{bipartite graph} $G = \{P, T, E\}$, where $P = \{1,2,...,N\}$ and $T = \{1,2,...,M\}$ are two independent and disjoint sets of nodes.
Further, we require that each edge in the \emph{edge set} $E = \{e_{ji}: i \in P, j \in T\}$ originates from $P$ and terminate at $T$.
As such, the considered network is \emph{directed}.

\begin{figure}[h!]
         \centering
         \includegraphics[scale=0.25]{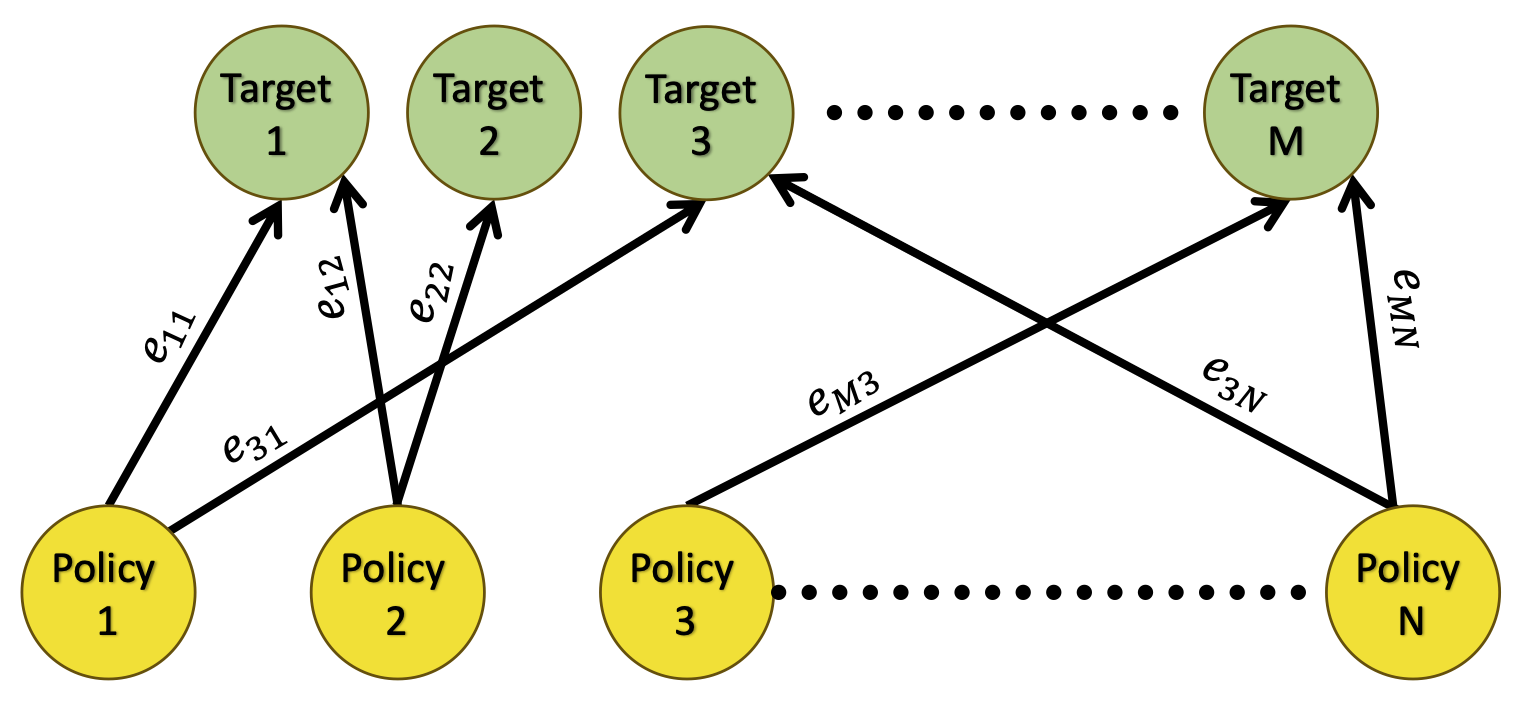}
    \caption{A schematic representation of a policy-target network. This two-mode network is directed, signed, and weighted.
    An edge $e_{ji}$ indicates that policy $i$ has a direct impact on target $j$.}
    \label{fig:bipartite_network}
\end{figure}

In our study, $P$ and  $T$ correspond, to the \emph{policy node set} and the \emph{target node set}, respectively.
The edge set $E$ represents the directed links that map policies to targets.
Consequently, the graph $G$ is designated as the \emph{policy-target network}; see    \cref{fig:bipartite_network}.
The outdegree of each policy $i \in P$ is denoted as $K_i^{out}$, while the indegree of each target $j \in T$ is denoted as $K_j^{in}$.

To encode the nature and strength of policy impacts on targets, we further augment the interlinkages with signed and weighted attributes.
More specifically, we assign a \emph{link coefficient} $c_{ji} \in [-1,1]$ to each existing edge $e_{ji}$ connecting policy $i$ to target $j$.
The sign of $c_{ji}$ (positive or negative) signifies whether the impact of policy $i$ on target $j$ is \emph{beneficial} or \emph{detrimental}, while the magnitude $|c_{ji}|$ quantifies the strength of this impact.
In summary, the policy-target network is a directed, signed, and weighted two-mode network, enabling a nuanced representation of complex interdependencies between policy actions and associated targets.

\subsection{Constructing the policy-target network using random distributions}
\label{sec:constructing_network}
Configuration models 
\citep{bollobas1980probabilistic} offer a set of statistical methods for constructing complex networks, providing a malleability in specifying the structure of connections among nodes.
In many real-world systems, the distribution of node degrees often adheres to a specific pattern, such as Gaussian, power-law or Poisson distributions.
Configuration models aim to capture and replicate such observed degree distributions in simulated networks.

The implementation of a configuration model typically involves two key steps: \emph{degree sequence generation} and \emph{edge assignment}.
To construct our policy-target network, we first generate an \emph{outdegree sequence} that specifies the number of \emph{outward edges} for each policy. 
Subsequently, we assign directed edges from policies to targets consistently with the outdegree sequence.
We assume that the specific target at which directed edges terminate is chosen \emph{uniformly at random}\footnote{However, one can modify our approach to accommodate distributions of the indegree sequence.}.
Therefore, for large networks, \emph{indegree sequence} for targets is expected to follow a normal distribution for large $M$.
Finally, since we are dealing with a weighted network, we also add a third step in which a link coefficient is assigned for each existing edge.

We carefully choose appropriate random distributions to generate the outdegrees of policies (a discrete random variable) and link coefficients for the policy-target interlinkages (a continuous random variable).
Furthermore, we will ensure that all generated random variables satisfy the following two criteria:
\begin{enumerate}
    \item \emph{independent and identically distributed}~(i.i.d), and
    \item within bounded intervals.
\end{enumerate}
To be flexible in our approach of constructing our network, we adopt \emph{generalized} random distributions capable of capturing a broad spectrum of distributional shapes, including symmetric, peaked and heavy-tailed profiles.
To maintain such flexibility and adhere to the two criteria outlined above, we employ modified Beta and Beta-binomial distributions to sample continuous and discrete random variables, respectively.

\subsubsection*{Continuous Beta distribution}
Consider a continuous random variable $X$ sampled from a \emph{Beta probability density function}:
\begin{equation*}
    X  \sim \ \text{Beta}(x;\alpha, \beta) = \frac{x^{\alpha-1} \cdot (1-x)^{\beta-1}}{\mathcal{B}(\alpha, \beta)}, 
\end{equation*}
where $\mathcal{B}(.)$ is the \emph{Beta function}, and $\alpha >0$ and $\beta >0$ are two shape parameters.
The random variable $X$ is supported by the continuous interval $[0,1]$.

Now consider a rescaled random variable $Y$ supported by the continuous interval $[y_{min},y_{max}]$.
We can obtain this variable by performing the following linear transformation
\begin{equation*}
    Y = (y_{max}-y_{min}) X + y_{min}.
\end{equation*}
The expected value of the rescaled variable $Y$ is given by
\begin{equation}
\label{eq:expected_v}
    \mu = E[Y] = \dfrac{(y_{max}-y_{min})\alpha}{\alpha+\beta} + y_{min},
\end{equation}
where $y_{min}$ and $y_{max}$ are the upper and lower bounds of $Y$ respectively.
Rearranging this equation allows us to express $\alpha$ as a function of $\beta$, $\mu$ and $y_{min}$ and $y_{max}$:
\begin{equation}
\label{eq:g_alpha}
    \alpha = \widetilde{\alpha}(\beta,\mu,y_{min}, y_{max}) = \dfrac{\beta(\mu-y_{min})}{y_{max}-\mu},
\end{equation}
where the average parameter $\mu$ must be chosen so that $y_{min} < \mu < y_{max}$.

Therefore, we can generate a bounded continuous random variable $Y$ using the following \emph{modified Beta distribution}:
\begin{equation}
\label{eq:final_continuous}
    Y = (y_{max}-y_{min}) X + y_{min}, \text{ where }  X \sim  \ \text{Beta}\bigl(x;\widetilde{\alpha}(\beta,\mu,y_{min},y_{max}), \beta \bigr),
\end{equation}
The random variable $Y$ can now be finely tuned by controlling four \emph{prescribed parameters}, namely, $\beta, \ \mu, \ y_{min},$ and $y_{max}$.
These prescribed parameters represent the scale parameter, the desired average, the lower bound and upper bound of the variable, respectively.

\subsubsection*{Discrete Beta-binomial distribution}
Consider a discrete random variable $X$ sampled from the \emph{Beta-binomial probability mass function}:
\begin{equation*}
    X  \sim \ \text{BetaBin}(x;\alpha, \beta, n) = \binom n x \ \frac{\mathcal{B}(x + \alpha, n - x+ \beta)}{\mathcal{B}(\alpha, \beta)}, 
\end{equation*}
where $\mathcal{B}(.)$ is the Beta function, $\alpha >0$ and $\beta >0$ are shape parameters, and $n \in \mathbb{N}$ represents the limit of the discrete support interval, i.e., $x \in \{0,1,...,n\}$.

Now consider a rescaled random variable $Y$ supported by the discrete interval $\{y_{min}, y_{min} + 1,...,y_{max}\}$.
To obtain $Y$, we first set 
\begin{equation*}
    n = y_{max} - y_{min}.
\end{equation*}
Then, we perform the following linear transformation
\begin{equation*}
    Y = X + y_{min}.
\end{equation*}
The expected value $\mu$ of the rescaled discrete variable $Y$ is given by equation    \cref{eq:expected_v}.
Therefore, we can express the scale parameter $\alpha$ in terms of $\beta$, $\mu$ and $y_{min}$ and $y_{max}$ using formula \cref{eq:g_alpha}, just as for the continuous case.

Hence, we can generate a bounded discrete random variable $Y$ using the following \emph{modified Beta-binomial distribution}:
\begin{equation}
\label{eq:final_discrete}
    Y =  X + y_{min}, \text{ where }  X \sim  \ \text{BetaBin}\bigl(x;\widetilde{\alpha}(\beta,\mu,y_{min},y_{max}), \beta, y_{max} - y_{min} \bigr),
\end{equation}
As for the continuous case, the random discrete variable $Y$ can now be finely tuned by controlling four prescribed parameters $\beta, \ \mu, \ y_{min},$ and $y_{max}$.
These parameters represent the scale parameter, the desired average, the lower bound and upper bound of the variable, respectively.

\subsubsection*{Policy outdegrees}
In the realm of policymaking, it is commonplace to observe that many policies suffer from inadequate planning, resulting in unintended consequences and an unforeseen rise in the number of targets influenced by each policy. 
Additionally, certain policies may naturally emerge as central hubs with higher outdegrees, exerting greater influence on the overall policy landscape. 
For these reasons, it is important to account for the diversity of policies in terms of their number of influenced targets.
Consequently, there is necessity for a non-uniform distribution of policies' outdegrees within our policy-target network.
To model the influence of policies, we configure a broad spectrum of statistical distributions of policy outdegrees by exploring various parameter combinations of the modified Beta-binomial distribution.

We begin by defining the outdegree of each policy $i$, denoted with the parameter $K^{out}_i \in \{1,2,...,M\}$, as the number of targets influenced by this policy.
The discrete bounded interval $\{1,2,...,M\}$ is based on the assumption that each policy must influence at least $1$ and at most $M$ targets. 
Then, we define the following \emph{outdegree sequence}: 
\begin{equation}
\label{eq:outdegree_array}
    \mathbf{K}^{out} = [K_1^{out},K_2^{out},...,K_N^{out}],
\end{equation}
which comprises the outdegrees of policies $1,2,...,N$.
To configure an outdegree sequence $\mathbf{K}^{out}$, we adopt the modified Beta-binomial distribution    \cref{eq:final_discrete}.
Specifically, we sample the outdegrees from the following distribution
\begin{equation}
\label{eq:outdegree_dist}
    Y_k =  X_k + 1, \text{ where }  X_k \sim  \ \text{BetaBin}\bigl(x;\widetilde{\alpha}(\beta_k,\mu_k,1,M), \beta_k, M - 1 \bigr),
\end{equation}
The discrete random variable $Y_k$ is determined by the three prescribed parameters $\beta_k$, $\mu_k$, and $M$, representing the scale parameter, the desired average of policy outdegrees, and number of targets, respectively.
By fixing two of the three parameters and varying the remaining parameter, one can generate a one-parameter family of Beta-binomial distributions parameterized by the varied parameter.
Additionally, we can estimate the \emph{network density} by taking the ratio of prescribed average outdegree and the number of targets: $\rho=\mu_k/M$.

\begin{figure}[h!]
     \centering
     \begin{subfigure}[b]{0.5\textwidth}
         \centering
         \includegraphics[scale=1.0]{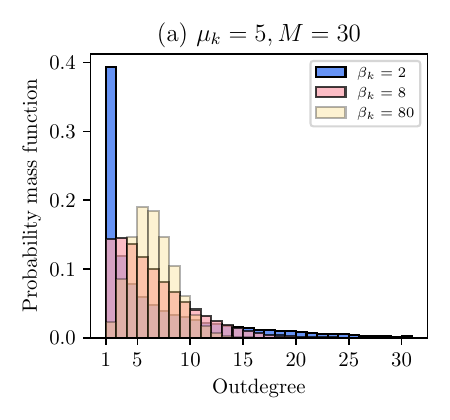}
     \end{subfigure}
     \hfill
     \begin{subfigure}[b]{0.45\textwidth}
     \centering
         \includegraphics[scale=1.0]{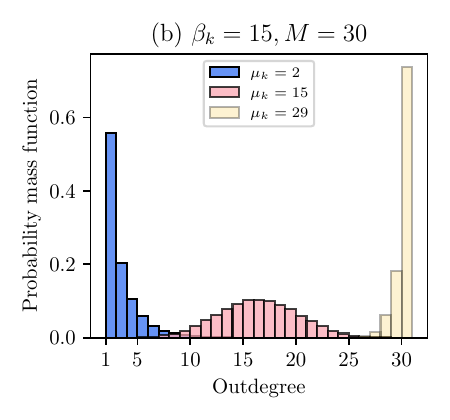}
     \end{subfigure}
      \caption{An illustration of two one-parameter families of outdegree distributions generated from equation    \cref{eq:outdegree_dist}.
      Panel~(a) shows a family of outdegree distributions with $M=30$, $\mu_k=5$ and varying values of $\beta_k$.
      Low, intermediate and high values of $\beta_k$, respectively, give rise to power-law-like, Poisson-like and near-symmetric distributions.
      Panel (b) exhibits a family of outdegree distributions with $M=30$, $\beta_k=15$ and varying values of $\mu_k$.
      For this value of $\beta_k$. low, intermediate and high values of $\mu_k$, respectively, result in to positively-skewed, near-symmetric, and negatively-skewed distributions.}
    \label{fig:outdegree_beta_distributions}
\end{figure}

 \Cref{fig:outdegree_beta_distributions} showcases two examples of one-parameter families of outdegree distributions generated from the discrete rescaled random variable $Y_k$ defined by equation    \cref{eq:outdegree_dist} with $M=30$ targets.
In \cref{fig:outdegree_beta_distributions}(a) where the prescribed average outdegree $\mu_k = 5$ is fixed, the $\beta_k$-dependent family of profiles range from a power-law-like distribution ($\beta_k = 2$), to a Poisson-like ($\beta_k = 8$) and a near-symmetric distributions ($\beta_k = 80$).
In panel (b) where the scale parameter $\beta_k = 15$ is fixed, the $\mu_k$-dependent profiles are near symmetric for intermediate  $\mu_k$ values (e.g. $\mu_k = 15$) but can also vary from positively-skewed (e.g., $\mu_k=2$) to negatively-skewed (e.g., $\mu_k=29$) for very small and large $\mu_k$ values, respectively.
Each of the two shown one-parameter family of configurations allows for a swift comparison along a continuous spectrum of network topologies explored by varying a single parameter.

\subsubsection*{Link coefficients and interaction matrix}

For every directed link from policy $i$ to target $j$, we assign a \emph{link coefficient} $c_{ji} \in [-1,1]$ using a continuous random distribution.
The sign and magnitude of this value represent the nature and intensity of the underlying link, respectively.
As such, the continuous interval $[-1,1]$ represents a range from the most detrimental ($-1$) to most beneficial ($+1$) impacts of policies on targets.
To this end, we construct following $M$-by-$N$ \emph{interaction matrix}:
\begin{align}  
\label{eq:interaction_matrix}
  \mathbf{C} =
  \left[ {\begin{array}{cccc}
    c_{11} & c_{12} & \cdots & c_{1N}\\
    c_{21} & c_{22} & \cdots & c_{2N}\\
    \vdots & \vdots & \ddots & \vdots\\
    c_{M1} & c_{M2} & \cdots & c_{MN}\\
  \end{array} } \right],
\end{align}  
where each link coefficient $c_{ji}$ quantifies the impact of policy $i$ on target $j$.
In the absence of a connection from policy $i$ to target $j$, the corresponding coefficient $c_{ji}$ is set to zero.
We note that, since this matrix serves as an algebraic representation of the policy-target network, the density of this matrix is equivalent to the network's density $\rho = \mu_k/M$.

We also use $R_j \in [-1,1]^N$ to denote the $j$th row of $\mathbf{C}$.
Hence, we can also express the interaction matrix as follows:
\begin{equation}
\label{eq:matrix_rows}
    \mathbf{C} = [R_1,R_2,...,R_M]^T.
\end{equation}
Alternatively, the interaction matrix can be expressed in terms of its columns:
\begin{equation}
\label{eq:matrix_columns}
    \mathbf{C} = [C_1,C_2,...,C_N].
\end{equation}
Each column $C_i \in [-1,1]^M$ captures the overall impact of policy $i$ on the influenced targets,  i.e., the number and magnitude of positive coefficients relative to the negative coefficients in the $i$th column of $\mathbf{C}$.

To configure link coefficients, we adopt the modified Beta distribution    \cref{eq:final_continuous}.
Specifically, we sample each link coefficient from the following distribution
\begin{equation}
\label{eq:signed_coeff}
    Y_c = 2 X_c - 1, \text{ where }  X_c \sim  \ \text{Beta}\bigl(x;\widetilde{\alpha}(\beta_c,\mu_c,-1,1), \beta_c \bigr).
\end{equation}
The constant quantities $\beta_c$ and $\mu_c$ represent the distribution's scale and average parameters, respectively. 
We also refer to the parameter $\mu_c$ as the \emph{overall efficacy} of policies, accounting not only for the number of beneficial and detrimental interlinkages but also their weights.
Note that the continuous random variable $Y_c$ is supported by the bounded continuous interval $[-1,1]$. 

\begin{figure}[h!]
     \centering
     \begin{subfigure}[b]{0.5\textwidth}
         \centering
         \includegraphics[scale=1]{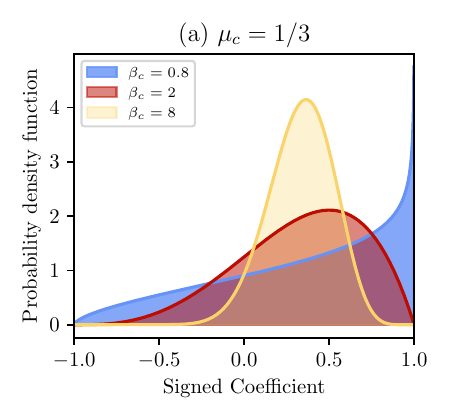}
     \end{subfigure}
          \hfill
     \begin{subfigure}[b]{0.45\textwidth}
     \centering
         \includegraphics[scale=1]{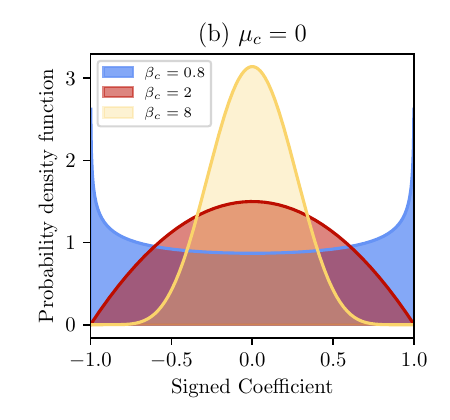}
     \end{subfigure}
      \caption{An illustration of two one-parameter families of link coefficient distributions generated from equation    \cref{eq:signed_coeff}, with (a) $\mu_c = 1/3$ and (b) $\mu_c=0$.
      Varying the underpinning parameters $\mu_c$ and $\beta_c$ enables flexible choices from a wide range of distributions including normal, heavy-tailed, Poisson, and convex distributions.}
    \label{fig:coeff_dist}
\end{figure}

   \cref{fig:coeff_dist} displays a range of illustrative examples of the Beta probability distributions for sampling link coefficients.
Panel (a) visualizes the case $\mu_c = 1/3$, where interlinkages are assumed to be predominantly beneficial.
The shape of the probability density function can be further modulated by fine-tuning the scale parameter $\beta_c$, giving rise to a spectrum of profiles ranging from negatively skewed to symmetric.
In contrast, panel (b) depicts the case where $\mu_c = 0$, which corresponds to an equilibrium between beneficial and detrimental interlinkages.
In this case, adjusting $\beta_c$ results in altering the concavity of the symmetric profiles.

\subsubsection*{Constructing the network}
To this end, we construct our policy-target network by taking the following steps.
First, we sample the outdegree sequence $\mathbf{K}^{out}= [K_1^{out},K_2^{out},...,K_N^{out}]$ using the corresponding random distribution    \cref{eq:outdegree_dist}.
Second, we assign each policy $i$ to $K_i^{out}$ distinct targets, selected uniformly at random.
Third, for each existing directed link $e_{ji}$ connecting a policy $i$ to a target $j$, we assign a link coefficient $c_{ji} \in [-1,1]$ sampled from the corresponding random distribution    \cref{eq:signed_coeff}.
In the absence of a directed link from policy $i$ to target $j$, the corresponding coefficient $c_{ji}$ is set to zero.

\begin{figure}[h!]
     \centering
     \begin{subfigure}[b]{0.13\textwidth}
         \includegraphics[scale=0.7]{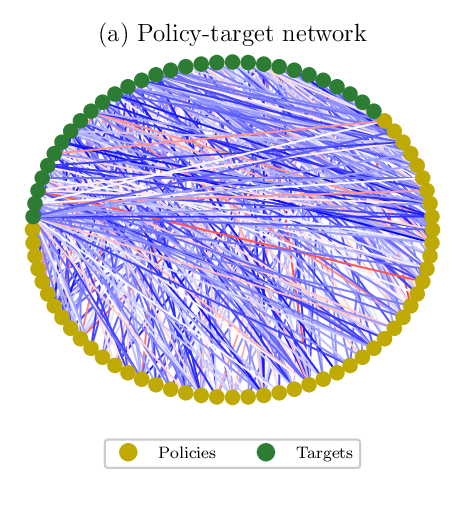}
     \end{subfigure}
          \hfill
     \begin{subfigure}[b]{0.65\textwidth}
         \includegraphics[scale=0.84]{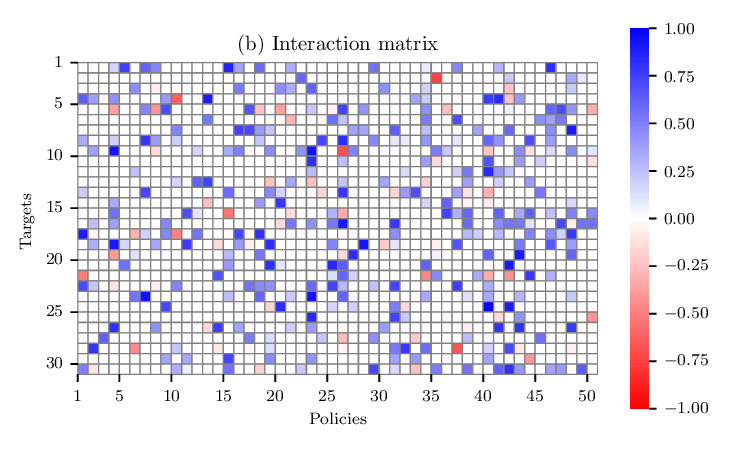}
     \end{subfigure}
      \caption{A demonstration of (a) a policy-target network and (b) a heatmap representation of the corresponding interaction matrix $\mathbf{C}$, constructed with $N=50$ policies and $M=30$ targets.
      This system was generated using an outdegree distribution    \cref{eq:final_discrete} with $\mu_k=7$ and $\beta_k =15$, and a link coefficient distribution    \cref{eq:signed_coeff} with $\mu_c=1/3$ and $\beta_c =2$.
      Colorbar indicates the value of the each link coefficient $c_{ji}$.}
    \label{fig:network_showcase}
\end{figure}

The result of this three-step procedure is a directed, signed and weighted two-mode network that links policies to targets.
   \cref{fig:network_showcase} showcases a policy-target network with $N=50$ policies and $M=30$ targets.  
In panel (a), the network's yellow and green nodes represent policies and targets, respectively. 
The directed links are color-coded according to the value of their respective link coefficient; see the color bar on the right for reference.

   \cref{fig:network_showcase} (b) displays a heatmap visualization of the corresponding
$M$-by-$N$ interaction matrix.
Each cell $(j,i)$ represents the direct influence of policy $i$ on target $j$, quantified by the link coefficient $c_{ji}$.
For example, consider policy 1, represented by the first column from the left. 
This policy has an outdegree of $K_1^{out} = 7$, indicating that it influences seven distinct targets.
Further, policy $1$ exerts a beneficial impact on targets 4, 9, 13, 17, 22, and 30, negative ramifications on target 21, and no influence on the remaining targets.

To modify the network density $\rho$ one can vary the average outdegree parameter $\mu_k$.
Meanwhile, the skewness of the outdegree distribution can be modified by fine-tuning the underpinning shape parameter $\beta_k$.
Similarly, the distribution of the link coefficients can be customized by adjusting the parameters $\mu_c$ and $\beta_c$.

It is important to emphasize that fixing the values of input parameters does not uniquely determine the realization of the interaction matrix $\mathbf{C}$. 
Due to the stochastic nature of the sampling process, repeated draws from the same underlying statistical distributions will generally produce different matrix realizations.
While the resulting $\mathbf{C}$ realizations of the same statistical distributions may differ quantitatively in their specific entries and link configurations, they are expected to remain qualitatively consistent in terms of their structural and statistical properties.
Nonetheless, as the system size increases, random realizations of $\mathbf{C}$  are expected to converge to the prescribed statistical distributions.

\subsection{Policy arrays and resource constraints}
The deployment of policies normally requires a meticulous plan for resourcing and budgeting.
In turn, policymakers are tasked with attaining objectives within the confines of a predetermined budget.
Therefore, we operate under two fundamental assumptions: (a) each policy incurs a specific cost, and (b) there exists a zero-sum competition among policies for a finite pool of public funds.
These assumptions underline the necessity for a judicious resource allocation, where the optimality of one policy may entail the reallocation from less workable policies.
Our goal is to numerically simulate a range of temporal scenarios that represent various sequences of interventions, emulating the iterative nature of policymaking and the consequential impacts on social and environmental targets.

In the classical NK fitness model, each `genome' is represented by an ordered sequence of `binary genes' that are either turned off or on. 
In a simple analogy to the context of our policy landscape, binary policies could either be activated or deactivated. 
However, to fine tune the impact of each policy, we define a non-binary alphabet parameter $A$. 
As such, the minimum and maximum possible share of resources assigned for each policy are given by $0$ and $A-1$, respectively.
For a balance between computational efficiency and the fine tuning of resource allocation, we need to choose an intermediate value of $A$.
Specifically, we fix $A=5$ throughout this work.

To represent the resource \emph{allocation} for each policy $i$, we define the discrete variable 
\begin{equation}
\label{eq:policy_allocation}
    x_i \in \{0,1,...,A-1\},
\end{equation}
capped by the upper ceiling $A-1$.
This variable is a dimensionless abstraction that represents a wide variety of resources, including labor, capital, material \citep{velez2023material}, energy \citep{grubler2018low}, land use, and technology, as well as negative utilities such as pollution and carbon emissions \citep{jackson2009prosperity}.
We note that there is a significant practical challenge in standardizing this metric across different areas; see, for example, Table 1 in \citep{schlesier2024measuring}.

Further, we define the \emph{policy array} (or \emph{policy configuration}) as
\begin{equation}
    \mathbf{X} = [x_1, x_2, \ldots, x_N],
\end{equation}
which represents the distribution of resources among policies.
In the context of multi-objective optimization, an instantiation of policy array is referred to as a \emph{solution} or \emph{decision} \citep{zhou2011multiobjective}.
When a policy array is traced over time, it is denoted as $\mathbf{X}[t]$ where `$t$' represents the adaptive step or simply `time'. 
Indeed, we view the policy array as a \emph{dynamic variable}
However, we omit the argument $[t]$ for notational convenience throughout most of the document.

Without any resource constraints, we operate within the \emph{full policy space}:
$$\mathcal{S} = \{0,1,...,A-1\}^N.$$
In optimization literature, $\mathcal{S}$ is often referred to as the \emph{solution space}.
Note that the size of the full policy space, given by $A^N$, grows exponentially with $N$.

To account for resource limitations as well as the inherent environmental constraints, we first define the following \emph{cost function} (also known as \emph{loss function}) of a policy array ${\mathbf{X}}$:
\begin{equation*}
    L({\mathbf{X}}) =
    \sum_{i=1}^{N} x_i,
    \label{eq:cost_function}
\end{equation*}
which represents the total amount of resources necessary for implementing the array's constituent policies.
Second, we assume that the \emph{total budget constraint}, denoted as $B_T$, is always fixed.
It follows that the \emph{average per-policy budget} $B = B_T/N$ also holds a constant value.
Hence, $B_T$ (or equivalently $B$) will be treated as a key input parameter that serves as a restriction to the amount of available resources \citep{desing2020ecological, desing2020circular}.
Mathematically, this restriction is imposed by the following inequality:
\begin{equation}
    L({\mathbf{X}}) \leq B_T \quad \text{or equivalently} \quad L({\mathbf{X}}) \leq B/N.
    \label{eq:constraint}
\end{equation}
We also use the term \emph{constraint boundary} to refer to the hypersurface $L({\mathbf{X}}) = B_T$.

To confine policy arrays within this mathematical constraint, we need to operate within the \emph{feasible region} $\mathcal{S}_B \subset \mathcal{S}$, defined as the set of all policy arrays satisfying the inequality above.
More formally, we define the feasible region as the following:
\begin{equation}
   \mathcal{S}_B :=  \{\mathbf{X}: \mathbf{X} \in \mathcal{S} \quad \text{and} \quad L({\mathbf{X}}) \leq B_T \}.
   \label{eq:feasible_region}
\end{equation}

\subsection{Relative importance of targets}
It is crucial to recognize that not all targets carry equal significance.
This is why we take a more refined approach that incorporates the \emph{relative importance} of targets\footnote{
It is also important to acknowledge that prioritization ranking is inherently subjective and context-sensitive.
In practice, adopting a participatory approach is indispensable.
This approach entails an inclusive engagement with all stakeholders and empowers local communities with agency in decision-making.
Furthermore, the prioritization of targets should never be static.
Regular review and adaptation are essential to maintain relevance and effectiveness in attaining the desired goals.
}.
Namely, we adopt a \emph{prioritization ranking} approach \citep{huan2022multi, tremblay2021systemic} to rank targets based on their importance and priority.
Each target $j$ is assigned a \emph{relative importance weight} $w_j \in \{1,2,...,10\}$, where 1 and 10 represent lowest and highest priority, respectively. 
We define the corresponding \emph{relative importance array}:
\begin{equation*}
    \mathbf{w} = [w_1,w_2,...,w_M]^T.
\end{equation*}
To configure a relative importance array $\mathbf{w}$, we employ the modified Beta-binomial distribution    \cref{eq:final_discrete}.
Specifically, each weight is sampled from the following distribution:
\begin{equation}
\label{eq:relative_importance}
    Y_w =  X_w + 1, \text{ where }  X_w \sim  \ \text{BetaBin}\bigl(x;\widetilde{\alpha}(\beta_w,\mu_w,1,10), \beta_w,  9 \bigr).
\end{equation}
To fine tune this random distribution, adjustments can be made to $\beta_w$ and $\mu_w$, representing the shape and average parameters of the random variable $Y_w$, respectively. 
Additionally, for a given relative importance array $\mathbf{w}$, we define an associated \emph{weight-normalizing constant}:
\begin{equation}
    \sigma(\mathbf{w}) = \sum_{j=1}^{M} w_j,
\end{equation}
which we will utilize in during the optimization step.

\subsection{Individual performance functions}

We characterize each target $j$ by an \emph{individual performance function} $f_j \in [0,1]$.
In the literature of fitness landscape models, such function has been referred to as the \emph{fitness contribution function} \citep{bull2022nonbinary, kauffman1989nk} or \emph{payoff function} \citep{baumann2024network}.
In the context of multi-objective optimization, an individual performance function is known as an \emph{objective function}.

We calculate the individual performance function for each target $j$ using the formula
\begin{equation}
    f_j(\mathbf{X}) =  \dfrac{1}{2} + \dfrac{1}{2} \mathcal{H}(\mathbf{X}, \eta) = \dfrac{1}{2} + \dfrac{1}{2}\tanh \Big( \dfrac{\eta}{(A-1) K^{in}_j} \ R_j \cdot \mathbf{X} \Big), 
\label{eq:contribution_nonlinear}
\end{equation}
evaluated at the policy array $\mathbf{X}$, where $R_j$ is the $j$th row of the interaction matrix $\mathbf{C}$ as expressed in~   \cref{eq:matrix_rows}.
The constant quantity $K^{in}_j$ reflects the indegree of target $j$, determined by the non-zero entries in $R_j$.
The
nonlinear function $\mathcal{H}(\mathbf{X}, \eta)$ is an \emph{activation function} of a logistic type, in which we set the scaling factor $\eta$ to 3.
This nonlinear scaling ensures that all individual performance functions remain within the continuous interval $[0,1]$. 
The interpretation of this approach is that the performance of target $j$ improves when resources are allocated to beneficial policies and deteriorates when allocated to detrimental ones.

In addition, we define the \emph{performance array} 
\begin{equation}
    \mathbf{f}(\mathbf{X}, \mathbf{C}) = \Bigl[f_1 (\mathbf{X}), f_2(\mathbf{X}), \ldots, f_M(\mathbf{X})\Bigr],
\end{equation}
as the array of all individual performance functions, evaluated at a given policy array $\mathbf{X}$.

\subsection{Multi-objective optimization}
When comparing two policy arrays $\mathbf{X}^{(1)}$ and $\mathbf{X}^{(2)}$, one needs to compare between individual components within their respective performance arrays $\mathbf{f}(\mathbf{X}^{(1)}, \mathbf{C})$ and $\mathbf{f}(\mathbf{X}^{(2)}, \mathbf{C})$. 
This comparison entails assessing performance across a multitude of competing targets and evaluating the underlying tradeoffs.
The process of identifying solutions that optimally balance the performance of multiple competing objectives is called \emph{multi-objective optimization}.

Typically, multi-objective optimization methods rely on the concept of \emph{Pareto optimality}. 
The core principle of Pareto optimality is that improving one or more objectives is not permitted if it comes at the expense of hindering another objective.
However, as the number of competing objectives increases, achieving Pareto dominance becomes increasingly difficult. 
More precisely, when optimising a large number of objectives, most Pareto comparisons end up in a tie and ultimately leading to a \emph{deadlock} where none of the compared arrays is \emph{Pareto dominant}.
Consequently, the time evolution of a Pareto local search is often obstructed by the high propensity of ending in a deadlock.
Additionally, Pareto optimality treats all objectives equally and does not accommodate for hierarchical importance of targets

To circumvent these challenges, we adopt an alternative approach, namely, \emph{scalarization} (or \emph{aggregation}), which is a widely-used method for multi-objective optimization \citep{gunantara2018review, giagkiozis2015methods,marler2010weighted}.
In this approach, the overall state of the system is evaluated via a scalarized function that incorporates the relative importance of targets.
Specifically, we define the \emph{overall performance} as a weighted sum of all individual performance functions:
\begin{equation}
    F(\mathbf{X}) = \dfrac{1}{\sigma(\mathbf{w})} \sum_{j=1}^{M} w_j  f_j(\mathbf{X}, R_j),
    \label{eq:matrix_form}
\end{equation}
evaluated at a given policy array $\mathbf{X}$.
The weight-normalizing constant $\sigma(\mathbf{w})$ ensures that $F(\mathbf{X}) \in [0,1]$.

\subsection{Optimization algorithm} \label{sec:opt_algorithm}
To optimize the overall performance $F(\mathbf{X})$ within the constraint of the resource budget $B_T$, we implement a combinatorial optimization algorithm that operates iteratively to seek a set of optimal solutions $\mathbf{X}$.
At each adaptive step $t$, our algorithm explores a local neighborhood within the feasible region, seeking improvements in the overall performance performance.
We opt for a local hill-climbing method in lieu of global search methods for several compelling reasons:
\begin{itemize} 
    \item The local search method enables the characterization of the performance landscape, by tracing the approximate number of local peaks (local optima) and capturing their heights and locations \citep{nowak2015analysis}.
    \item It reanimates the incremental nature of the \emph{social learning process} \citep{baumann2024network}, also referred to as \emph{exploration process} or \emph{innovation}, which typically occurs within a local search space.
    This local approach contrasts with the idealized assumption of omniscient decision-makers capable of globally surveying the entire policy space.  
    \item It emulates the temporally evolving nature of policy interventions, where adjustments are typically made through iterative refinements \citep{fontanari2015exploring}. 
    \item It facilities parallel evaluation of multiple intervention trajectories originating from distinct initial conditions, akin to principles found in particle swarm optimization of social interactions \citep{kennedy1999minds}.
        \item 
    \cite{weinberger1996np} has proved that the original NK fitness landscape model is
    \emph{NP-complete}, i.e., the task of finding global optima computationally intractable for large $N$ .
    Since we add several complexities to the original NK model (e.g., non-binary policy strings, two-mode network and herirarchical importance of objectives), we conjecture that our optimization problem is also likely to be NP-complete.
\end{itemize}

To enhance reproducibility and analytical consistency, we employ a \emph{deterministic} algorithm for local search \citep{kauffman1987towards, nowak2015analysis}.
However, while the optimization path within each simulation is deterministic, stochastic elements are introduced via the random initialization of both the interaction matrix $\mathbf{C}$ and the initial policy configuration $\mathbf{X}[0]$ at the beginning of each run. 
Thus, to obtain robust and generalizable insights, we perform an ensemble of a sufficiently large number independent simulations (e.g., $n=100$).
This ensemble-based approach facilitates a comprehensive exploration of variability across system configurations.


\subsubsection*{The algorithm}
We now formalize our local search algorithm, which predicates on navigating local neighborhoods of the performance landscape.
Firstly, we refer to two policy arrays $\mathbf{X}^{(1)}$ and $\mathbf{X}^{(2)}$ as \emph{adjacent neighbors} if they differ in a single component by exactly one unit.
Since all policy arrays are composed of non-negative integers, a policy array $\mathbf{X}^{(1)}$ is an adjacent neighbor of $\mathbf{X}^{(2)}$ if and only if their Euclidean distance $d(\cdot)$ is equal to one, that is, 
\begin{equation*}
    d(\mathbf{X}^{(1)} - \mathbf{X}^{(2)}) = 1.
\end{equation*}

Secondly, consider a dynamic policy array $\mathbf{X}[t]$ at some adaptive step $t$.
We define an
\emph{adjacent neighborhood} 
$\mathcal{N}(\mathbf{X}[t])$ of the current policy array $\mathbf{X}[t]$ as the set of all its neighbors, including the state itself.
Formally:
\begin{equation*}
    \mathcal{N}(\mathbf{X}[t]) = \Bigl\{ \mathbf{X}: \mathbf{X} \in \mathcal{S} \text{ and } d(\mathbf{X}[t] - \mathbf{X}) = 0 \text{ or } 1 \Bigr\}.
\end{equation*}

Thirdly, when navigating a neighborhood, we need to operate within the confines of the feasible region $\mathcal{S}_B$ defined in equation    \cref{eq:feasible_region}.
Thus, we introduce the concept of the \emph{feasible neighborhood}, denoted as $\mathcal{N}_B(\mathbf{X}[t])$, which represents the set intersection of $\mathcal{N}(\mathbf{X}[t])$ and $\mathcal{S}_B$, and can be expressed as follows:
\begin{equation*}
    \mathcal{N}_B(\mathbf{X}[t]) = 
    \Bigl\{\mathbf{X}:\mathbf{X} \in \mathcal{N}(\mathbf{X}[t]) \quad \text{and} \quad L({\mathbf{X}}) \leq B_T \Bigr\}.
\end{equation*}

Finally, an initial condition $\mathbf{X}[0] \in \mathcal{S}_B$ of the system at time $t=0$ is chosen at random. 
At each subsequent time step $t+1$, the \emph{adaptive process} (also termed the \emph{selection process}) depends on the set of neighbors of the current state $\mathbf{X}[t]$.
Specifically, we seek the policy array with the highest overall performance $F(\mathbf{X}, \mathbf{w})$ within the current feasible neighborhood $\mathcal{N}_B(\mathbf{X}[t])$:
\begin{equation}
    \mathbf{X}[t+1] = \argmax_{\mathbf{X} \in \mathcal{N}_B(\mathbf{X}[t])} F(\mathbf{X},\mathbf{w}), \quad \mathbf{X}[0] \in \mathcal{S}_B.
    \label{eq:transitionUndecimal}
\end{equation}
This adaptive process continues until a \emph{local optimum} is reached \citep{nowak2015analysis}.
We define a local optimum as the policy array that outperforms all its neighbors within $\mathcal{S}_B$.
In contrast, a \emph{global optimum} is defined as the policy array with highest performance in the entire feasible region $\mathcal{S}_B$.
A transition from policy array to its neighbor is referred to as an \emph{intervention}---in analogy to the term mutation in the context of fitness landscapes.

We use the term \emph{trajectory} (or \emph{pathway}) to describe a temporal sequence of interventions.
The following is an example of a trajectory that starts from some initial condition $\mathbf{X}[0]\in \mathcal{S}_B$, proceeds through a number of interventions, before settling at some local optimum $\mathbf{X}[T]$:
\begin{equation}
    \mathbf{\Gamma} = [\mathbf{X}[0], \mathbf{X}[1],..., \mathbf{X}[T]],
\end{equation}
where $T$ denotes the number of adaptive steps taken to reach a local optimum, also referred to as \emph{convergence time} or \emph{trajectory length}.
Note that the convergence time of a trajectory depends on its initial condition.
\footnote{In dynamical systems terminology, our time-evolution system is autonomous, i.e., the adaptive step does not explicitly depend on the value of $t$.
Furthermore, our \emph{initial value problem}    \cref{eq:transitionUndecimal} guarantees one and only one solution (a unique solution) for any initial condition $\mathbf{X}[0] \in \mathcal{S}_B$.
We also note that our adaptive function in    \cref{eq:transitionUndecimal} is a non-surjective non-injective function $g:\mathcal{S}_B \to \mathcal{S}_B$, that is, not all elements of the co-domain is associated with an argument and some images could be associated with more than one argument.}


\subsection{Rationale of selecting parameter values}

To calibrate estimate baseline values of the system parameters, we first give a brief overview of the relevant studies from the literature.
\cite{muneeb2022assessing} developed a multicriteria decision-making model for optimizing energy resource allocation for economic and environmental goals in India.
\cite{weitz2018towards} investigated a weighted target-target matrix to assess interactions among 34 SDG targets in Sweden.
Similarly, \cite{mccollum2018connecting} analyzed how energy-related SDG targets are related other targets.
\cite{van2019analysing} created an SDG interaction matrix using expert surveys and text mining. The UN Global Sustainable Development Report 2019 \citep{messerli2019global} also developed a framework to construct a 17-by-17 matrix for SDG interactions.
This framework was also detailed in a peer-reviewed article authored by \cite{pham2020interactions}.
In addition, \cite{breu2021begin} downscaled this framework to the national level with a case study focusing on Switzerland.

Furthermore, a major framework for mapping the interlinkages of all 169 SDG targets in 27 countries was developed in \citep{moinuddin2017sustainable}.
\cite{dawes2020sustainable} introduced a linear dynamical system to model the progress of SDGs, while \cite{algunaibet2019powering} studied sustainable energy systems through planetary boundaries.
\cite{kaper2019modeling} surveyed various modeling frameworks for food systems within the less-studied framework of Doughnut Economics.
\cite{kroll2019sustainable} examined trade-offs and synergies between SDGs globally, and \cite{lusseau2019income} introduced the concept of the ``sustainome" to map SDG interactions using a vast World Bank dataset in conjunction with network analysis.
\cite{moallemi2020global} explored medium- and long-term pathways for achieving SDGs through scenario assessments of a systems-dynamics model.
On the policy front, \cite{lowder2023food} presented and analyzed a large-scale database of international policy levers that address food systems and \cite{stechemesser2024climate} estimated the impact of 1,500 carbon-emission policies across 41 countries and highlighted that successful policies are often part of a coordinated mix rather than standalone interventions.

Across these studies, we have identified the following overarching patterns.
First, co-benefits consistently outweigh trade-offs, indicating that most interlinkages within the system are synergistic in nature.
In the context of our policy-target network, the signed link coefficients will be inferred to be predominantly positive in order to reflect the pervasiveness of co-beneficial interactions.
Second, the outdegree distributions in the studied systems generally follow a power-law pattern (this pattern is also called the Pareto principle).
These distributions were also typically characterized by a relatively small average degree.
This observation supports our hypothesis that a few high-degree policies tend to exert significant influence, while the majority of policies are expected to have limited impact.
This concentration of influence suggests that a scale-free network is a suitable choice for our analytic framework.

Additionally, the number of policies can be treated as a more flexible parameter, as the influence of each policy depends on the scale of its constituent interventions.
Similarly, the number of targets is contingent on their scope, potentially ranging from 21 to 169, corresponding to the number of dimensions in the Doughnut Economics framework \citep{raworth2017doughnut} and the total number of SDG targets \citep{SDGs2015}, respectively.

A crucial oversight in previous studies has been the omission of targets' relative importance.
However, a recent study \citep{kaufmann2024democratic} examined how citizens and policymakers rank \emph{urban policy issues} (e.g., cost of living, public health, and biodiversity). 
Despite the significant observed discrepancies between the two perspectives in this study, both sets of rankings followed an exponential distribution, indicating a hierarchical prioritization of desired outcomes.
Our model reflects this hierarchical nature of targets' prioritization to ensure that resources and efforts are allocated efficiently.

Building on these insights, we examine the behavior of our model by exploring an ensemble of various statistical scenarios, calibrated to align with existing findings.
To reflect the dominance of co-benefits over tradeoffs, we focus on the case of a positive policy efficacy ($\mu_c>0$). 
This imposed condition indicates a net-positive suite of policy interventions.  
To introduce asymmetry into the distribution of policy impacts, we adopt a negatively-skewed distribution of link coefficients by setting $\beta_c=2$ as a baseline value.

To model a sparse and scale-free policy-target network, we select relatively small values for $\mu_k$ and $\beta_k$.
For a hierarchical ranking of targets, we choose $\mu_w = 8$ and $\beta=15$ as baseline values.
For the baseline number of policies and targets, we set $N=100$ and $M=30$, respectively. 
Unless explicitly stated otherwise, all simulations and analyses presented in this study utilize the baseline parameter values summarized in    \cref{table:parameters}.
Finally, it is important to treat these parameter values with flexibility to accommodate any uncertainties.

\section{Results} \label{sec:Results}
In this section, we begin by showcasing our optimization algorithm by presenting a prototype performance landscape.
We then investigate the system's behavior under various budgetary scenarios.
Finally, we perform a simulation-based sensitivity analysis to assess the robustness of our model and identify the key parameters and driving factors of the system.

\begin{figure}[h!]
         \centering
         \includegraphics[scale=0.85]{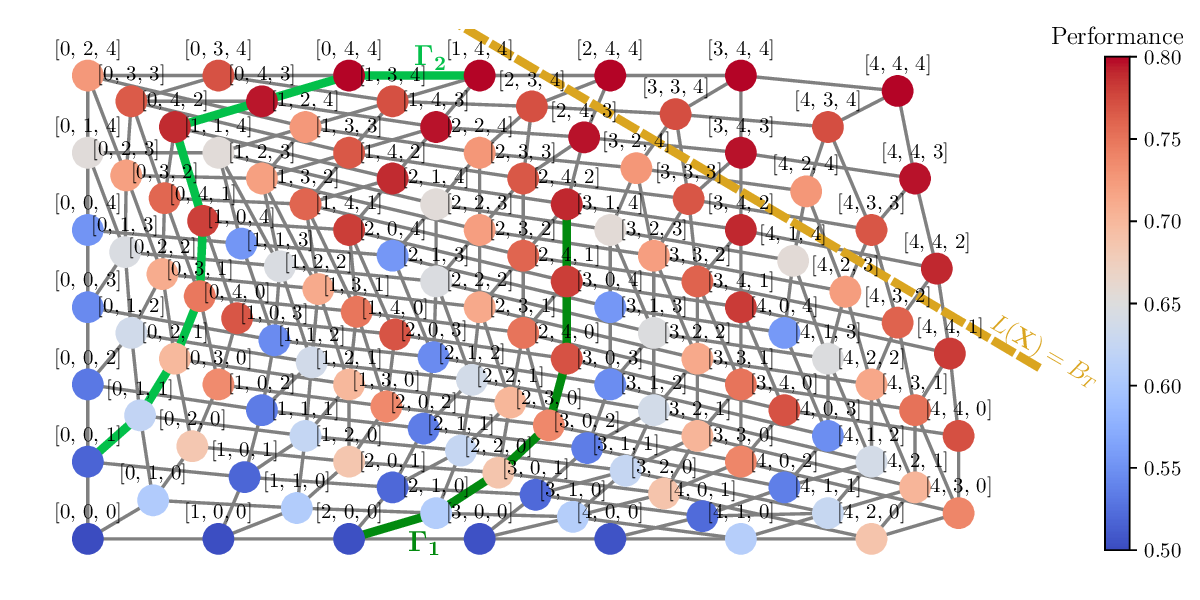}
    \caption{A schematic illustration of performance landscape with $N=3$ policies, $M=30$ targets, alphabet of $A=5$, and a budget constraint of $B_T=9.5$.
    Each policy array $\mathbf{X} = [x_1,x_2,x_3]$ is visualized as a vertex and colored according to its corresponding performance.
    Gray edges linking two policy arrays indicate that they are adjacent neighbors. 
    The diagonal line $L(\mathbf{X}) = B_T $ represents the constraint boundary of the feasible region $\mathcal{S}_B$;  $\mathcal{S}_B$ comprises the subregion below this line and the line itself.
    The first trajectory of interventions, denoted $\mathbf{\Gamma_1}$, starts from a low-performance state and terminates at a local optimum with higher performance. 
    The second trajectory $\mathbf{\Gamma_2}$ starts from a low-performance state and reaches the global optimum within $\mathcal{S}_B$.
    }
\label{fig:performance_landscape}
\end{figure}

\subsection{Navigating the performance landscape}
To showcase our model and algorithm, we visualize the full policy space $\mathcal{S}$ as a performance landscape.
In classical NK fitness landscapes  \citep{kauffman1989nk}, a genotype is typically visualized as a single locus in a two-dimensional continuous domain, and the fitness of each genotype is represented as an elevation.
In this representation, local maxima and minima are visualized as peaks and valleys, respectively.
In the context of our framework, we aim to visualize each policy array as a single vertex in a two-dimensional discrete performance landscape, with its corresponding overall performance function encoded by a topographic colormap; see also  \citep{de2014empirical}.

A schematic visualization of a performance landscape with $N=3$ policies, $M=30$ targets, alphabet of $A=5$, and a budget of $B_T=9.5$, is illustrated in \Cref{fig:performance_landscape}.
Each policy array $\mathbf{X} = [x_1,x_2,x_3]$ represents the amount of resources allocated for policies 1, 2 and 3, reflecting the respective scales of implementation.
The aim is to simultaneously optimize the performance of $M=30$ targets (not shown) within a finite budget of $B_T=9.5$ (orange dashed line).
Each variable $x_i$, corresponding to the resource allocation in policy $i$, is capped by the prescribed upper ceiling $A-1=4$, and so must be within the discrete set $\{0,1,2,3,4\}$.
The performance $F(\mathbf{X})$ of the system for each policy array $\mathbf{X}$ is indicated by the color bar, signifying the `elevation' of the corresponding array within the landscape.
A gray link between two arrays represents their adjacency as neighbors.

The full policy space $\mathcal{S}$ is represented by the entire two-dimensional plane.
The cost function $L(\mathbf{X})$ of a policy array $\mathbf{X}$ quantifies its total resource expenditure. 
Policy arrays that invoke
the same cost function are aligned along the same \emph{negative diagonal} (i.e., diagonal line with a negative slope).
As such, each negative diagonal, serves as a \emph{contour line} of a fixed cost function.
The diagonal line $L(\mathbf{X}) = B_T$ defines the constraint boundary of the feasible region $\mathcal{S}_B$ for the choice of $B_T=9.5$.
The budget parameter $B_T$ abstracts how the system is constrained by finite amount of available resources and the constraints of the ecological boundaries  \citep{rockstrom2009planetary}.  
The feasible region $\mathcal{S}_B$ is represented by the subregion below this line, including the line itself. 
In other words, the aggregate amount of resources for each policy array below or on this line does not exceed $B_T=9.5$.
Note that one would obtain the same result for any choice $B_T\in[9,10)$.

Two examples of optimization trajectories evolving from two distinct initial conditions are shown in \cref{fig:performance_landscape}.
The first trajectory $\mathbf{\Gamma_1}$ (dark green) originates from the policy array $[2,0,0]$ and encounters a sequence of interventions before terminating at the policy array $[2,4,3]$.
During every adaptive step, resources are reallocated between policies in pursuit of a higher performing state within the neighborhood of the current state.
This (greedy) hill climbing process is repeated until a local optimum, namely, $[2,4,3]$, is reached  \citep{nowak2015analysis}.
The second trajectory $\mathbf{\Gamma_2}$ (light green) starts from another initial condition, namely, $[0,0,1]$, and navigates the landscape before converging to the local optimum $[1,4,4]$, which also happens to be the global optimum of the feasible region $\mathcal{S}_B$. 
Both trajectories are confined within $\mathcal{S}_B$. 
By removing the budget constraint, i.e., setting $B_T=12$, both trajectories can evolve further in time and are expected to reach the global optimum $[4,4,4]$ of the full policy space $\mathcal{S}$.

By varying the system parameters, one can modify the size and topology of the performance landscape as well as the elevations within the landscape.
For example, increasing the parameters $N$ and $A$ expands the size of the landscape since the total number of policy arrays is given by $A^N$. 
Changing the number of targets $M$ and their relative importance or modifying the statistical distributions of the underlying policy-target network may have a significant effect on the behavior of the system.
These factors influence the number of peaks and valleys, and consequently determine the topography of the landscape.
Larger systems are more challenging to visualize, which is why we chose a small $N$ for the demonstration of the performance landscape shown in \Cref{fig:performance_landscape}.

\begin{figure}[h!]
     \centering
     \begin{subfigure}[b]{0.32\textwidth}
         \centering
         \includegraphics[scale=0.7]{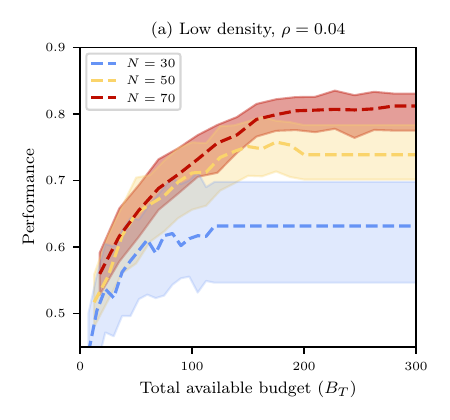}
     \end{subfigure}
     \hfill
     \begin{subfigure}[b]{0.32\textwidth}
     \centering
    \includegraphics[scale=0.7]{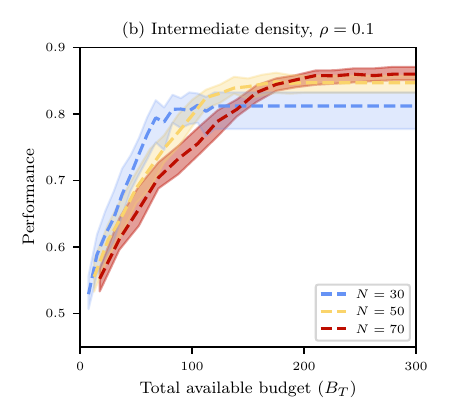}
     \end{subfigure}
     \hfill
     \begin{subfigure}[b]{0.32\textwidth}
     \centering
    \includegraphics[scale=0.7]{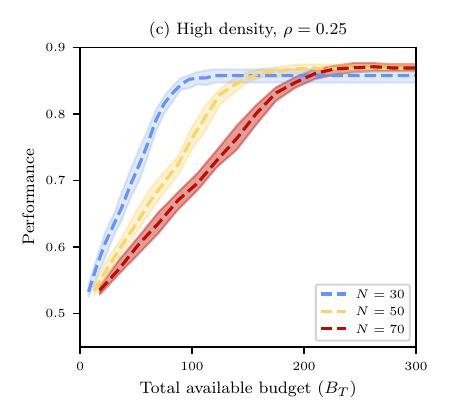}
     \end{subfigure}
    \caption{Performance trends for increasing total budget $B_T$ across different policy scenarios ($N=30$, 50, and 70) with varying levels of network density: (a) low, $\rho=0.04$, (b) intermediate, $\rho=0.1$ and (c) high, $\rho=0.25$.
    All other parameters are given by the baseline values outlined in the \cref{table:parameters}, except for the number of targets, fixed at $M=50$.
    The lower and upper limits of each trend represent the 25th and 75th percentiles, respectively, with the median shown as a dashed line.
    For every $B_T$ value, an ensemble of $n=100$ simulations were conducted, each corresponding to random realization (newly generated matrix $\mathbf{C}$).
    At each simulation, starting from a randomly generated initial condition, the system is evolved until reaching a local optimum, as described in \cref{sec:opt_algorithm}.
    }
    \label{fig:budget_stats}
\end{figure}

\subsection{System performance under various scenarios}
The primary focus of our analysis is the evaluation of the overall performance of the system subject to changes in the amount of disposable resources or total budget, denoted as $B_T$.
We also aim to examine how the system responds to variations in both the number of policies represented by $N$ and the policy-target network density $\rho = \mu_k/M$.

We start by examining the trends of the overall performance as the budget parameter $B_T$ is tuned.
\Cref{fig:budget_stats} demonstrates three policy scenarios with $N=30$, $50$ and $70$, each simulated under three different network density levels: low, intermediate and high.
For each scenario, 100 experimental simulations were conducted.
Across all scenarios, we observe a general trend where overall performance consistently improves as $B_T$ is increased.
Nevertheless, as $B_T$ reaches large values, the performance begins to level off, and the marginal gain (i.e., rate of incremental improvement) gradually declines.
This phenomenon is widely known as a \emph{point of diminishing returns} (PDR)  \citep{knight1944diminishing, Riera-Prunera2014}, where the marginal gains in the output plateau or even decline as incremental units of input are added.

This observation aligns with earlier studies  \citep{raworth2017doughnut,jackson2009prosperity,stern2007economics,kallis2025post,meadows1972limits}. Grounded in empirical evidence, these studies have demonstrated that perpetual economic growth ultimately 
leads to stagnation in social and environmental outcomes.
While incremental growth in monetary investments may lead to improvements up to a certain PDR, exceeding this optimal PDR may result in inefficiencies and wasted resources.
Therefore, in a world with finite resources and an expanding global population, sustainable resourcing and budgeting practices are essential to preclude the overexploitation of resources, which could otherwise lead to system collapse or long-term adverse effects on both society and ecosystems  \citep{kallis2025post, meadows1972limits}.



To this end, we evaluate the system's behavior subject to variations in the number of policies $N$.
Three key findings emerge from our analysis, as depicted in \cref{fig:budget_stats}.
First, for small values of $N$, the system stagnates at relatively low budget $B_T$.
For instance, the blue scenario ($N=30$) in \cref{fig:budget_stats}(a) reaches a PDR at $B_T=120$ compared to the yellow scenario ($N=50$) at $B_T=200$.
Such premture stagnation of performance for small $N$ occurs because of the fixed ceiling on the resource allocations per policy, set at $A-1=4$ (see equation \cref{eq:policy_allocation}).
Consequently, the overall performance cannot improve beyond the ceiling $B \geq A-1$, or equivalently $B_T \geq N(A-1)$.
For example, when $N=30$, it is noticeable that increasing the $B_T$ beyond $N(A-1)=120$ yields no further improvement due to this ceiling.
Therefore, aligning the number of policies with the availability of resources is crucial to avoid inefficiencies.

Second, expanding the number of policies $N$ is likely to enhance the system's robustness and mitigate uncertainties as it enables a scope for diversifying beneficial policies while simultaneously distributing the risks among multiple targets, rather than concentrating risks on a few outcomes.
Nonetheless, this process of widening the scope of policies must not come at the expense of the overall efficacy of policies---quantified by $\mu_c$.  
While it is theoretically viable to maintain a high value of $\mu_c$ as $N$ is increased, it may prove difficult in practice.

Third, increasing $N$ only pays off if the amount of allocated resources grows correspondingly.
We observe that the optimal number of policies $N$ depends on two main factors: (a) total amount of invested resources, symbolized by $B_T$, and (b) the density $\rho$ of the policy-target network---which is linearly scaled by the average outdegree $\mu_k$.
When the policy-target network is sparse (low $\rho$), increasing $N$ consistently improves the performance, as illustrated in \cref{fig:budget_stats}(a).
Conversely, in the case of a dense policy-target network (high $\rho$), large values of $N$ become less advantageous, especially for a low budget; see \cref{fig:budget_stats}(c).
For moderate values of $\rho$, the impact of $N$ becomes less consequential to the overall performance, as depicted in \cref{fig:budget_stats}(b).

In the classical NK landscape model, an optimal value of the interconnectivity parameter $K$ is typically found at an intermediate level  \citep{kauffman1989nk,bull2022nonbinary}.
In contrast, the findings of our modified landscape model suggests that the ideal value of the network density $\rho$ is contingent on the amount of disposable resources $B_T$ as well as the number of policies $N$.
With a high budget $B_T$, increasing $\rho$ generally yields better outcomes; on the other hand, a low $B_T$ necessitates a careful balance between $N$ and $\rho$.
Moreover, the variance of the recorded values of local optima is negatively correlated with the system's density $\rho$. 
In other words, the variance narrows for larger $\rho$.
It follows that a policy diversification strategy makes the system less sensitive to initial conditions and more robust to external shocks.

A notable caveat to this analysis is the empirical challenge of detecting unintended consequences, both negative and positive, which
often results in the underestimation of $\rho$.
To address this, we recommend the explicit distinction between advertent and inadvertent outcomes when mapping policies to targets.
Despite the practical challenge of the fine tuning of $\mu_c$ in empirical settings, our findings suggest that adopting holistic policies---those that are characterized with high $\rho$ and simultaneously tackle multiple targets---is likely to be instrumental and resource-efficient.
Nevertheless, the optimization process for such holistic policies may require an intricate analysis of intertwined tradeoffs, which can be cumbersome in practice.



\begin{figure}[h!]
     \centering
     \begin{subfigure}[b]{0.32\textwidth}
         \centering
         \includegraphics[scale=0.525]{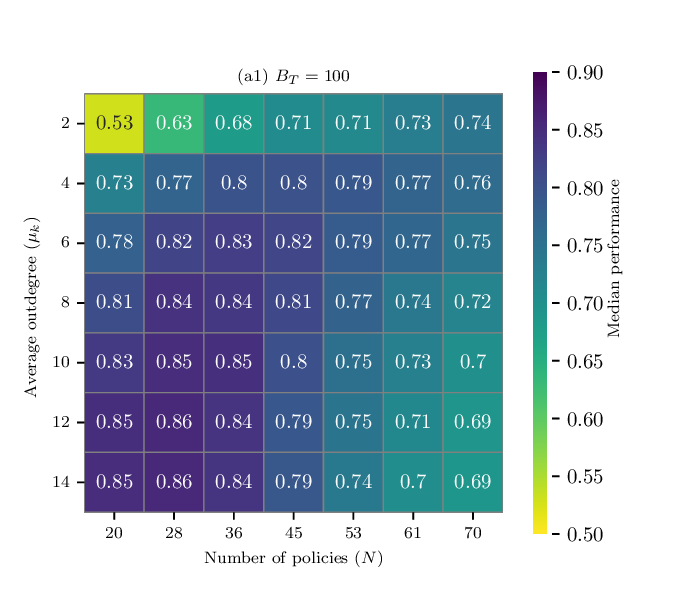}
     \end{subfigure}
     \hfill
     \begin{subfigure}[b]{0.32\textwidth}
     \centering
    \includegraphics[scale=0.525]{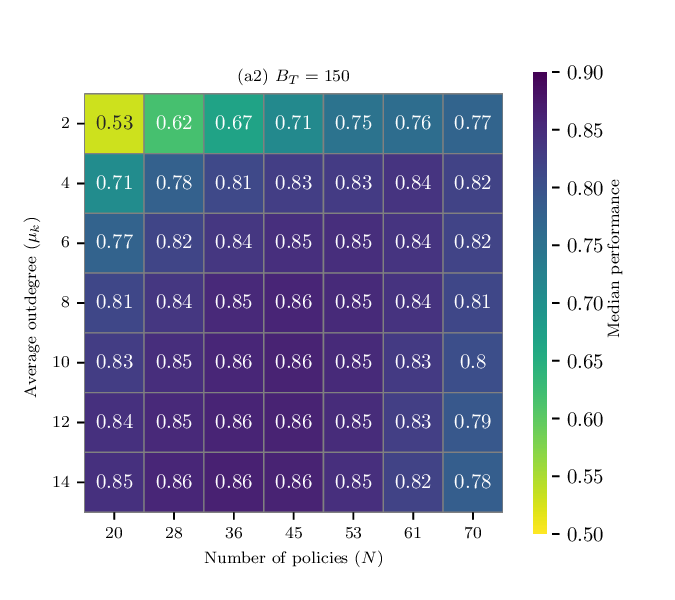}
     \end{subfigure}
     \hfill
     \begin{subfigure}[b]{0.32\textwidth}
     \centering
    \includegraphics[scale=0.525]{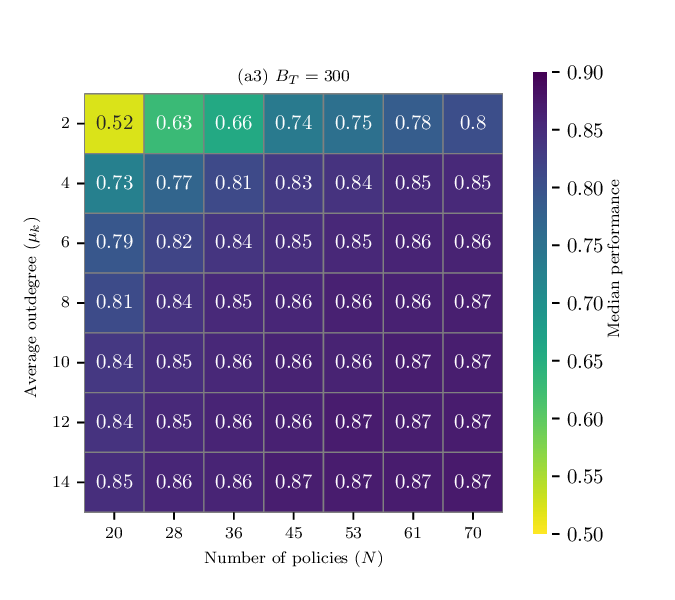}
     \end{subfigure}
     \begin{subfigure}[b]{0.23\textwidth}
     \centering
    \includegraphics[scale=0.525]{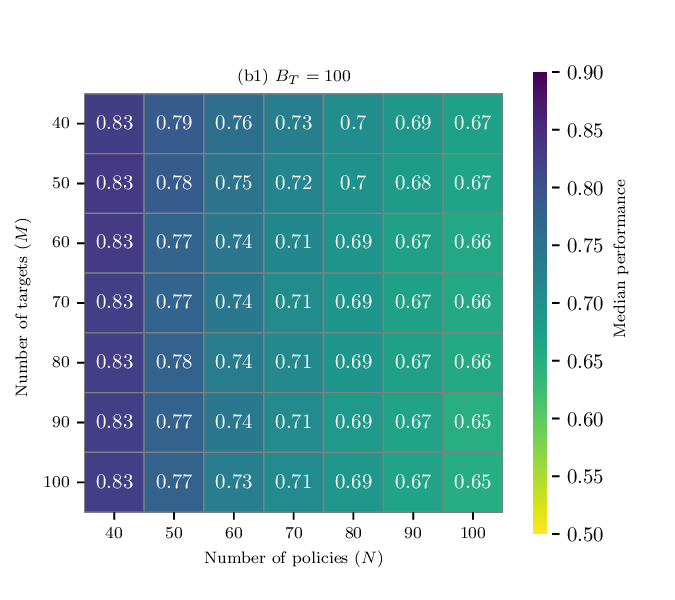}
     \end{subfigure}  
     \hfill
     \begin{subfigure}[b]{0.23\textwidth}
     \centering
    \includegraphics[scale=0.525]{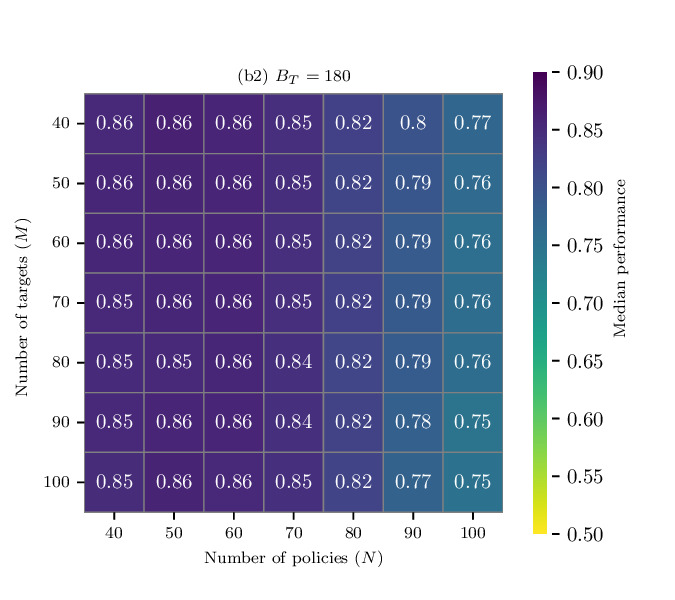}
     \end{subfigure}
     \hfill
     \begin{subfigure}[b]{0.34\textwidth}
     \centering
    \includegraphics[scale=0.525]{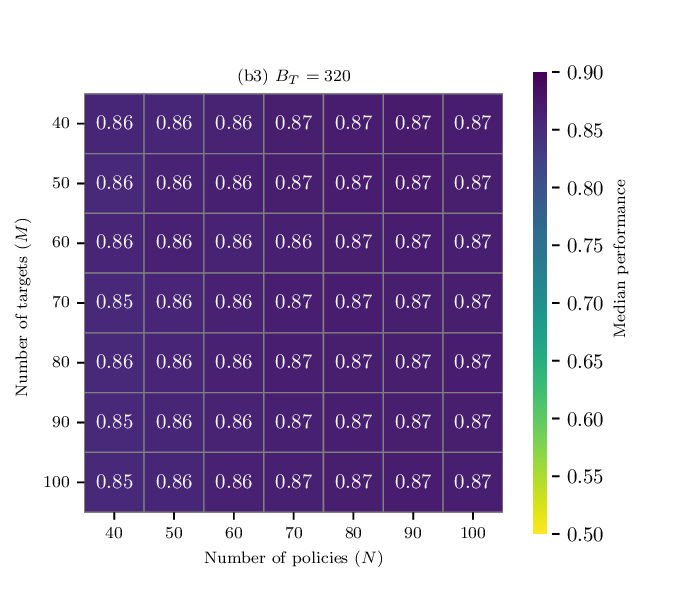}
     \end{subfigure}
     \begin{subfigure}[b]{0.23\textwidth}
     \centering
    \includegraphics[scale=0.515]{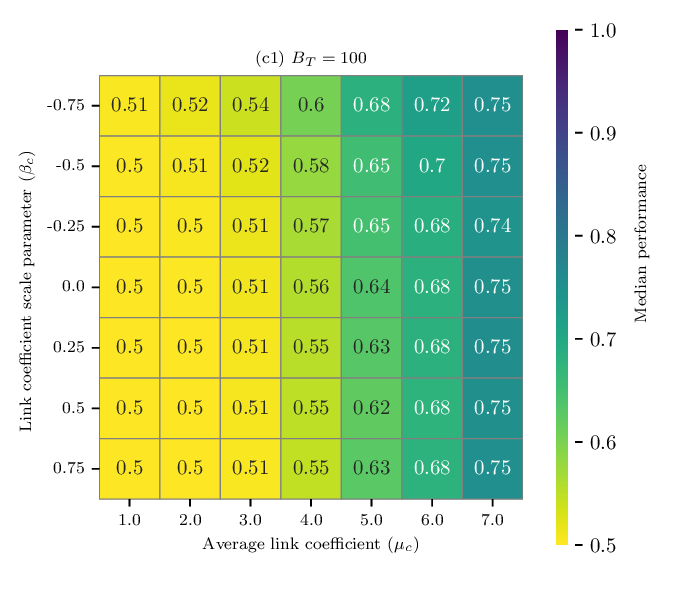}
     \end{subfigure}  
     \hfill
     \begin{subfigure}[b]{0.23\textwidth}
     \centering
    \includegraphics[scale=0.515]{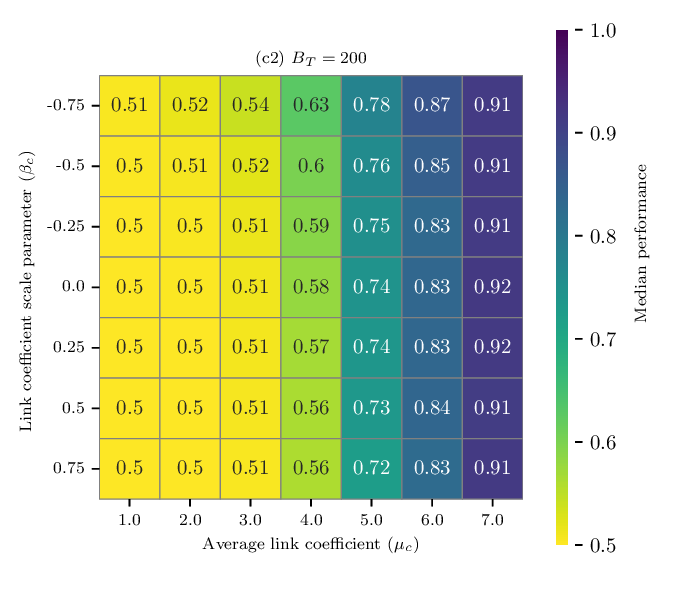}
     \end{subfigure}
     \hfill
     \begin{subfigure}[b]{0.34\textwidth}
     \centering
    \includegraphics[scale=0.515]{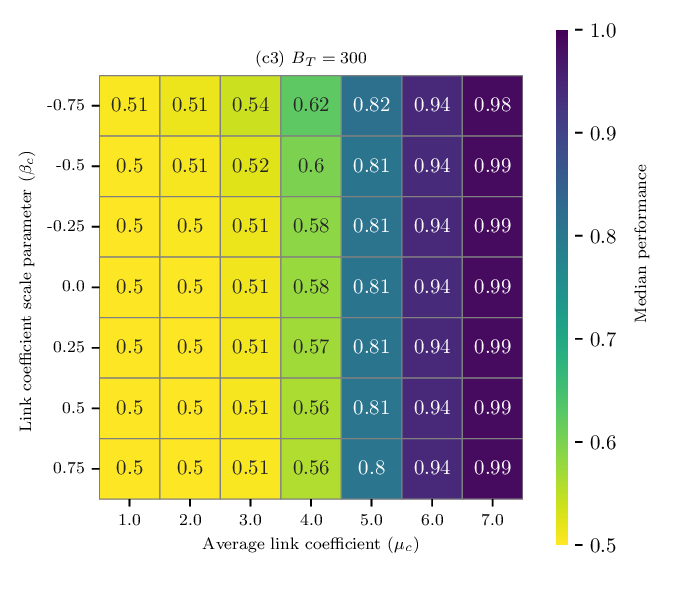}
     \end{subfigure}
    \caption{Sensitivity analysis for the the overall performance under varying input parameters.
    The heatmaps illustrate the interplay between (a) the number of policies $N$ and network density $\mu_k$, with $M=50$, (b) the number of policies $N$ and number of targets $M$, with fixed network density $\rho=\mu_k/M=1/6$, and (c) the average policy efficacy $\mu_c$ and the associated scale parameter $\beta_c$.
    The columns reveal the change of the system behavior as the total budget $B_T$ increases. 
    The value in each cell represents the median performance across 100 stochastic simulations.
    At each simulation, a new interaction matrix \textbf{C} and a random initial condition are first generated, then the system is evolved until reaching a local optimum, as described in \cref{sec:opt_algorithm}.
    Parameter values are as listed in \cref{table:parameters}, unless stated otherwise above.
    }
    \label{fig:sensitivity_analysis}
\end{figure}

\subsection{Sensitivity analysis}

To comprehensively assess the system's robustness and pinpoint the key factors in driving the system behavior, we conduct an extensive \emph{sensitivity analysis} across multiple system parameters.
This analysis aims to unravel how variations in different input parameters affect the system's overall performance.
We consider a discrete range of input parameter values,
which are summarized in \cref{table:parameters}.


Our sensitivity analysis begins by examining the optimal relationship between the number of policies $N$ and their average outdegrees $\mu_k$.
Panels (a1)--(a3) of \cref{fig:sensitivity_analysis} illustrate a transition of distinct patterns across three budgetary scenarios.
These results reveal that achieving optimal performance requires a calibrated balance between $N$ and $\mu_k$ as follows.

For a low budget (e.g., $B_T=100$), the highest-performance block (darkest in color) is positioned at the bottom left corner of the ($N,\mu_k$) parameter plane.
Under such tight resource constraint, the system benefits from a reduced number of policies, with the notable exception when $\mu_k$ is very low. 
For a very high budget (e.g., $B_T=300$), 
it is advantageous to boost the number of policies $N$, allowing for comprehensive utilization of resource abundance.
For this high-budget scenario, the optimal region in the ($N,\mu_k$) parameter plane shifts to the bottom right corner of the heatmap.
In the intermediate case of a moderate budget (e.g., $B_T=150$), the balance between $N$ and $\mu_k$ becomes more nuanced; the optimal region is concentrated around the center of the parameter space, suggesting that optimal outcomes require a meticulous consideration for deploying interventions.
It is important to reiterate that while decision-makers have direct control over $N$, the underlying network density (quantified by $\rho = \mu_k / M$) may inadvertently be underestimated due to the unforeseen emergence of unintended consequences of policies over targets.

Next, we investigate the relationship between the number of policies $N$ and number of targets $M$; see \cref{fig:sensitivity_analysis}(b).
To ensure consistency, network density is fixed with $\rho =1/6$.
The analysis suggests a significant system sensitivity to variations in $N$.
However, as discussed earlier, the sensitivity to $N$ diminishes as the total budget $B_T$ increases.
Remarkably, the system shows minimal sensitivity to changes in the number of targets $M$, given that the ratio $\rho = \mu_k/M$ is kept fixed. 
In other words, expanding the scale of desired outcomes is not expected to affect the system performance, as long as the average policy outdegree $\mu_k$ scales proportionally with the number of targets.
Indeed, when comparing \cref{fig:sensitivity_analysis}(a) with \cref{fig:sensitivity_analysis}(b), it appears that $\mu_k$ plays a more pivotal role than $M$ in attaining favorable outcomes. 
This suggests that policymakers need to focus on managing and maximizing holistic reach of policies when aspiring to a large number of objectives.

To this end, we examine the sensitivity of the system's behavior to the average policy efficacy $\mu_c$ and the associated scale parameter $\beta_c$. 
\Cref{fig:sensitivity_analysis}(c) displays a strong positive correlation between system performance and $\mu_c$.
In contrast, the results indicate a near-negligible sensitivity to variations in $\beta_c$.
The general pattern in the ($\mu_c$, $\beta_c$) parameter plane persists across different values of the budget $B_T$, indicating that this pattern is robust under various resource constraints.
This finding aligns with theoretical expectations that higher values of the efficacy parameter $\mu_c$ imply numerous and overwhelmingly beneficial policies, which in turn enhance a wide range of desired goals.

Finally, we have conducted preliminary analysis on the impact of the parameters $\mu_w$ and $\beta_w$, which together determine the random distribution of the relative importance weights of targets.
We found that neither parameter exhibited a statistically significant impact on the overall performance.

\section{Discussion} \label{sec:Discussion}

This manuscript introduces a generalizable analytical framework for monitoring and evaluating the impact of diverse policy scenarios on sustainability targets.
The complex relationships between policy levers and their potential outcomes were represented as a policy-target network.
Drawing on the seminal work on NK fitness landscapes  \citep{kauffman1989nk}, our framework takes the form of a multi-objective combinatorial optimization powered by an evolutionary algorithm and network analysis. 

The primary motivation behind this work is the recognition that progress toward sustainability often unfolds in silos.
For instance, environmental targets (e.g., reducing transport-based emissions) are frequently pursued independently from social initiatives (e.g., reducing inequalities).
Conversely, interventions that aim to advance socioeconomic outcomes may inadvertently undermine or exacerbate ecological vulnerabilities.  
Such compartmentalized planning can result in negative unintended consequences or untapped opportunities to leverage synergies between various goals.
To circumvent these challenges, a holistic approach  \citep{meadows2008thinking} is paramount.
Our framework provides a whole-system methodology for developing a decision-support toolkit for addressing problems
characterized by complex interdependencies.

\subsection{Methodological contributions}
The methodological framework proposed in this paper integrates methods from network science, statistical configurations, and evolutionary algorithms. 
This framework provides a systematic pipeline for a multi-objective optimization that goes beyond the traditional methods applied in NK fitness landscape models.


At the core of our framework is a two-mode, directed, weighted, and signed policy-target network.
This two-mode structure captures the complex relationships between interventions and intended outcomes---a notable departure from the one-mode networks typical of the classical NK model. 
Additionally, our approach can be readily extended to accommodate higher-order networks.

To construct the policy-target network, we employed configuration models  \citep{bollobas1980probabilistic}.
Specifically, we sampled the outdegrees of policies and link coefficients for the network using Beta (continuous) and Beta-binomial (discrete) distributions, respectively.
This approach accommodates a more flexible and systematic assignment of policy-target interlinkages, including their signs and magnitude. 
The policy-target network was also mathematically formalized via an interaction matrix. 


Furthermore, our framework introduces a method for allocating and distributing finite resources across policies via a non-binary alphabet cardinality  \citep{bull2022nonbinary, srivastava2023alphabet}.
This formulation allows for finer granularity in policy specification and resource deployment.
As such, the size of the full policy space ($\mathcal{S}$), given by $A^N$, expands at an exponential rate with the number of policies $N$. 
However, the size of the feasible region ($\mathcal{S}_B$) is capped by the limitations of the available resources and environmental constraints, parametrized by $B$. 
While the large size of the policy space may pose a computational challenge, our approach enables a fine-grained scaling of policies and interventions, moving beyond the limitations of traditional binary landscapes.

In classical NK models, the influence of input variables (e.g., genes) is often encoded using
randomly generated values  \citep{kauffman1989nk, kauffman1987towards, weinberger1996np, stoltzfus2006mutation}, while output indicators (e.g., fitness metrics) are typically treated with uniform importance.
In contrast, our approach offers a more transparent and structured evaluation of system performance by modeling both (a) the efficacy of each policy and (b) hierarchical ranking across competing targets.
The efficacy of policy decisions is captured through (individual) performance functions that compute the progress towards each target based on both the amount of resources allocated to relevant policies and the magnitude and direction of their influence.
The prioritization of targets is governed by relative-importance parameters.
These parameters were then combined to compute the aggregate system performance as a weighted average of individual performances, allowing a detailed prioritization of outcomes.

Our combinatorial optimization algorithm employs a deterministic, iterative hill-climbing approach  \citep{kauffman1987towards}.
This algorithm navigates the local neighborhood of a current policy configuration by making incremental adjustments, while respecting budget constraints.
This approach is powerful for visualizing and navigating complex policy landscapes (see \cref{fig:performance_landscape}).
Furthermore, the implemented algorithm supports the exploration of multiple potential intervention pathways by initializing from different policy configurations.
Unlike Pareto optimization, which often leads to decision deadlocks in high-dimensional objective spaces, our approach not only avoids such deadlocks but also incorporates hierarchical prioritization of targets to reconcile trade-offs more effectively.



\subsection{Main findings}
The main finding of our model suggests that increasing the amount of invested resources leads to a point of diminishing returns, as demonstrated in \cref{fig:budget_stats}. 
This observation aligns with established theoretical and empirical studies. 
Most notably, the seminal work by \cite{meadows1972limits} highlights the inherent limitations to perpetual growth, and theorizes the possibility of ``overshoot and collapse," a phenomenon that provides a mechanism for tipping points  \citep{mcbain2017long, ritchie2021overshooting}.
Furthermore, our results support the strong argument that \emph{prosperity}  \citep{jackson2009prosperity} \emph{thriving societies}  \citep{raworth2017doughnut} can be achieved without a direct reliance on indefinite growth. 

The presented analysis reveals that targets, even when directly symbiotic, can impede each other's progress due to the implicit tradeoffs embedded in policy-target interlinkages.
Such indirect competition is modeled in our framework through the way resources are allocated among counteracting policies. 
Our optimization model also suggests that the transient disruptions in the progress of a subset of targets are often necessary to accomplish long-term, system-wide, optimal outcomes.
While such delays and disruptions are necessary for sustainable transitions, they pose a practical challenge for governance where policymakers, pressured by immediate results, may abandon long-term strategies in favor of short-term gains that ultimately obstruct or undermine holistic progress.

The proposed framework was complemented with a thorough sensitivity analysis to identify the critical drivers of the system. 
The budget parameter $B_T$, emerged as the most influential factor.
This parameter not only symbolizes the amount of available resources, but also represents an abstraction of the constraints of environmental and planetary boundaries  \citep{rockstrom2009planetary}.
The second key determinant of the system is the overall policy efficacy $\mu_c$; the amplification of this parameter skews the distribution of policy efficacy, yielding higher performance outcomes.
The number of policies $N$ and the network density $\rho$ were found to be less pivotal compared to $B_T$ and $\mu_c$, though still relevant.
Notably, the overall efficacy of policies remains paramount. 
In practice, decision-makers may iteratively tune the number and scale of policies to strike an optimal balance between $N$ and $\mu_k$ as illustrated in Figure 3(a).
Furthermore, we found that the precise number of targets $M$ is of minimal significance when the ratio $\rho = \mu_k/M$ is kept constant. 
The parameters $\mu_w$ and $\beta_w$,  which govern the relative importance weights, exhibit negligible influence on the qualitative behavior of the system.

In real-world settings, policymakers have limited control over the number of outcomes influenced by the average policy---quantified by $\mu_k$.
This limitation arises from inherent blind spots and the potential for unpredictable repercussions.
Consequently, it is crucial to assess the net-sum impact of policies, by weighing the magnitude of synergistic impacts against any detrimental side effects, which is measured in our model the efficacy parameter $\mu_c$.
This warrants a holistic approach for policy design that seeks to maximize co-benefits while minimizing trade-offs.
This conclusion is verified by the findings of our simulations, where resources are promptly reallocated from antagonistic policies to more efficient alternatives.

\subsection{Caveats and challenges}
Although it has its strengths, our approach does not come without limitations and caveats.
First, while top-down policy interventions play a key role in driving urban transformation, they do not solely determine the system-wide performance.
Rather, outcomes are shaped by a collection of multiple interdependent factors, including citizen engagement, private sector involvement, and broader socioeconomic complex dynamics.
Currently, our team is developing a performance landscape model that integrates top-down and bottom-up approaches to examine the optimal degree of devolution in the decision-making process.

Another caveat stems from the economic asymmetries among nations in their capacity to pursue system transformations, indicating that our framework may not be universally applicable.
For instance, high-income countries and cities have access to more advanced technologies and established sustainable infrastructure enables multiple, cost-effective pathways to accelerate sustainable transitions.
Nonetheless, many of these pathways rely on the ongoing extraction of raw materials and rare earth elements, frequently sourced from less developed regions. 
This dependency raises pressing ethical concerns and necessitates fair trade standards and equitable practices throughout the global supply chain  \citep{nsude2024global}.


\subsection{Future recommendations}

Moving forward, a multidimensional approach to budgeting and resource allocation is essential.
This could be accomplished by disaggregating the budget into distinct categories---or `faucets'---of resources, such as labor, capital, land use, raw materials, machinery alongside resources that extend beyond exchange values (e.g., pollution and carbon emissions).
While our aggregate budget analysis provides insights into the impact of resource distribution on system performance, a more fine-tuned examination of how each of the intimately intertwined faucets of resources distinctly affect specific targets remains a a key area for future research.

While the proposed framework offers valuable guidance for decision making in the context of sustainable development, a reliable evaluation of the impact of policy interventions on sustainability targets presents a significant empirical challenge.
\cite{stechemesser2024climate} recently introduced a novel approach to evaluate climate policy measures seeking to reduce carbon emissions.
Although this study operated on the assumption that policy impacts can be captured only within a specified temporal window 
following implementation, the proposed approach nonetheless establishes a broadly applicable, data-driven method for policy evaluation. 
The integration of our analytical framework with this emerging approach for policy evaluation is a promising avenue for future research.




\section*{Author contributions}
CRH conceptualized the study, designed the methodology, wrote the programming code, performed simulations and analysis, and wrote the original draft of the manuscript. 
JC conceptualized the study, co-designed the methodology, reviewed the manuscript, and provided supervision.
All authors contributed to methodology, analysis, and reviewing the manuscript.

\section*{Funding and acknowledgment}
This study was funded by the Natural Environment Research Council as part of the Changing the Environment Programme (grant number NE/W005042/1).

\pagebreak

\appendix

\section{Glossary}\label{secA1}
\begin{itemize}
\item {\bf Target}: A specific objective that represents a desired social or environmental outcome. 
\item {\bf Policy}: A policy lever available to decision- and policy-makers in pursuit of sustainability targets.
\item {\bf Policy-target network}: A directed two-mode network that consists of two disjoint sets of nodes, representing $N$ policies and $M$ targets.
\item {\bf Interaction matrix} ($\mathbf{C}$):  An $M$-by-$N$ matrix representation of the policy-target network.
\item {\bf Link coefficient} ($c_{ji}$): A quantity that represents the nature (positive or negative) and the \\ magnitude (relative extent) of the impact of policy $i$ on target $j$.
\item {\bf Policy outdegree} ($K_i^{out}$):     Number of targets impacted by a policy $i$. 
\item {\bf Outdegree sequence} ($\mathbf{K}^{out}$):   An ordered array of policy outdegrees. 
\item {\bf Predominately synergistic (resp. antagonistic)}:  A policy-target interaction matrix where the vast majority of relationships are beneficial (resp. detrimental).
\item {\bf Policy allocation}  ($x_i$):  Number of discrete amount of resources allocated to a given policy $i$. 
\item {\bf Policy array}  ($\mathbf{X}$):  An ordered array of policy allocations---also called a `solution' in the context of multi-objective optimization.
\item {\bf Alphabet} ($A$):   The number of discrete amount (cardinality) of resources that can be allocated to policies.
\item {\bf Total budget} ($B_T$):  Total amount of disposable resources that can be distributed among policies. 
\item {\bf Per-policy budget} ($B$):  Average amount of disposable resources allocated per policy. i.e., $B = B_T/N$.  
\item {\bf Full policy space} ($\mathcal{S}$):  The set of all unique policy arrays. 
\item {\bf Feasible region} ($\mathcal{S}_B$): A subset of $\mathcal{S}$ that only includes policy arrays not exceeding the budget. 
\item {\bf Relative importance weight} ($w_j$):  A quantity that indicates the relative importance of a given target $j$. 
\item {\bf Relative importance  array} ($\mathbf{w}$):  An ordered array of relative importance weights. 
\item {\bf Individual performance function} ($f_j$): A measure of the performance (outcome) for a given target $j$---also called an `objective function' in the context of multi-objective optimization.
\item {\bf Performance array}  ($\mathbf{f}$)::   An ordered array of individual performance functions. 
\item {\bf Overall performance} ($F$): Aggregated sum of individual performances weighted according the relative importance of the impacted targets.
\item {\bf Adjacent neighbors}:  Two policy arrays with Euclidean distance of one. 
\item {\bf Adjacent neighborhood} ($\mathcal{N}(\mathbf{X})$): The set of neighbors of a given policy array $\mathbf{X}$. 
\item {\bf Feasible neighborhood} ($\mathcal{N}_B(\mathbf{X})$): The intersection set of $\mathcal{N}(\mathbf{X})$ and $\mathcal{S}_B$. 
\item {\bf Local optimum}: A policy array (`solutions') that performs better than all other arrays within its feasible neighborhood.
\item {\bf Global optimum}:  A policy array that performs better than all other arrays in the entire feasible region.

\end{itemize}

\pagebreak
\bibliographystyle{unsrt}
\bibliography{sn-bibliography}
\end{document}